\documentclass[twocolumn]{autart}

\usepackage{graphicx}
\usepackage{amsmath}
\usepackage{amssymb}
\usepackage{times,mathptm}
\usepackage{amsfonts}
\usepackage{enumerate}
\usepackage{longtable}
\usepackage{graphicx}
\usepackage{theorem}
\usepackage{multirow}
\usepackage{mdframed}
\usepackage[square,sort,comma,numbers]{natbib}
\usepackage{float}
\usepackage{color}

\newcommand{\CR}[1]{{#1}}
\DeclareMathOperator*{\argmin}{arg\,min}

\newtheorem{remark}{Remark} [section]
\newtheorem{lemma}{Lemma} [section]
\newtheorem{proposition}{Proposition} [section]
\newtheorem{assumption}{Assumption} [section]
\newtheorem{theorem}{Theorem} [section]

\begin{document}

\begin{frontmatter}
%\runtitle{Insert a suggested running title}  % Running title for regular 
                                              % papers but only if the title  
                                              % is over 5 words. Running title 
                                              % is not shown in output.

\title{An analysis of the SPARSEVA estimate for the finite sample data case \thanksref{footnoteinfo}}

\thanks[footnoteinfo]{The material in this paper was not presented at any conference.}

\author[First]{Huong Ha}\ead{huong.ha@uon.edu.au},
\author[First]{James S. Welsh}\ead{james.welsh@newcastle.edu.au},
\author[Second]{Cristian R. Rojas}\ead{cristian.rojas@ee.kth.se},
\author[Second]{Bo Wahlberg}\ead{bo.wahlberg@ee.kth.se}

\address[First]{School of Electrical Engineering and Computer Science, The University of Newcastle, Australia}
\address[Second]{Department of Automatic Control and ACCESS, School of Electrical Engineering, KTH Royal Institute of Technology, SE-100 44 Stockholm, Sweden}


\begin{abstract}
In this paper, we develop an upper bound for the SPARSEVA (SPARSe Estimation based on a VAlidation criterion) estimation error in a general scheme, i.e., when the cost function is strongly convex and the regularized norm is decomposable for a pair of subspaces. We show how this general bound can be applied to a sparse regression problem to obtain an upper bound for the traditional SPARSEVA problem. Numerical results are used to illustrate the effectiveness of the suggested bound.
\end{abstract}

\end{frontmatter}

\section{Introduction}

Regularization is a well known technique for estimating model parameters from measured input-output data. Its applications are in any fields that are related to constructing mathematical models from observed data, such as system identification, machine learning and econometrics. The idea of the regularization technique is to solve a convex optimization problem constructed from a cost function and a weighted regularizer (regularized M-estimators). There are various types of regularizers that have been suggested so far, such as the $l_1$ \cite{Tibshirani:96}, $l_2$ \cite{Tikhonov:77} and nuclear norms \cite{Fazel:01} \cite{Fazel:03}. 

During the last few decades, in the system identification community, regularization has been utilised extensively \cite{Pillonetto:14}\CR{, to impose properties of smoothness and sparsity in the estimated models (see, e.g., \cite{Ohlsson-10,Zorzi-Chiuso-17})}. Most of this work has focused on analysing the asymptotic properties of an estimator, i.e., when the length of the data goes to infinity. The purpose of this type of analysis is to evaluate the performance of the estimation method to determine if the estimate is acceptable. However, in practice, the data sample size for any estimation problem is always finite, hence, it is difficult to judge the performance of the estimated parameters based on asymptotic properties, especially when the data length is short.

Recently, a number of authors have published research (\cite{Bickel2009}, \cite{Candes2007}, \cite{Negahban2012}) aimed at analysing estimation error properties of the regularized M-estimators when the sample size of the data is finite. Specifically, they develop upper bounds on the estimation error for high dimensional problems, i.e., when the number of parameters is comparable to or larger than the sample size of the data. Most of these activities are from the statistics and machine learning communities. Among these works, the paper \cite{Negahban2012} provides a very elegant and interesting framework for establishing consistency and convergence rates of estimates obtained from a regularized procedure under high dimensional scaling. It determines a general upper bound for regularized M-estimators and then shows how it can be used to derive bounds for some specific scenarios.%However, this framework is only applicable for the case when the number of parameters is comparable, or larger, than the sample size of the data whereas in a typical system identification problem, the number of parameters is generally smaller than the sample size of the data.

Here in this paper we utilize the framework suggested in \cite{Negahban2012} to develop an upper bound for the estimation error of the M-estimators used in a system identification problem. Here, the M-estimator problems are implemented using the SPARSEVA (SPARSe Estimation based on a VAlidation criterion) framework \cite{Rojas:11}, \cite{Rojas:14}. \CR{The approach in \cite{Negahban2012} has been developed for penalized estimators, so it has to be suitably modified for SPARSEVA, which is not a penalized estimator, but the solution of a constrained optimization problem.} Our aim is to derive an upper bound for the estimation error of the general SPARSEVA estimate. We then apply this bound to a sparse linear regression problem to obtain an upper bound for the traditional SPARSEVA problem with some assumptions on the regression matrix. These assumptions can be considered as the price in order to derive the upper bound. In addition, we also provide numerical simulation results to illustrate the suggested bound of the SPARSEVA estimation error.

The paper is organized as follows. Section 2 formulates the problem. Section 3 provides definitions and properties required for the later analysis. The general bound for the SPARSEVA estimation error is then developed in Section 4. In Section 5, we apply the general bound to the special case when the model is cast in a linear regression framework. Section 6 illustrates the developed bound by numerical simulation. Finally, Section 7 provides conclusions.

\CR{
\subsection{Notation}

In this paper, we will use the following notation:
\begin{itemize}
\item
$f(x\vert 0,\sigma^2) = (2\pi\sigma^2)^{-1/2} \exp(-x^2/2\sigma^2)$ denotes the probability density function (pdf) of the Normal distribution  $\mathcal{N}(0,\sigma^2)$.
\item
$\chi^2_\beta(N)$ denotes the value that $P(X < \chi^2_\beta(N)) = 1-\beta$, where $X$ is Chi square distributed with $N$ degrees of freedom.
\end{itemize}
}

\section{Problem Formulation}

Let $Z_1^N = \{Z_1, ..., Z_N\} \CR{\in \mathcal{Z}^N}$ denote $N$ identically distributed observations with marginal distribution $\mathbb{P}$ \CR{in $\mathcal{Z} \subseteq \mathbb{R}^k$}. $\mathcal{L}: \mathbb{R}^n\times\mathcal{Z}^N \rightarrow \mathbb{R}$ denotes a convex and differentiable cost function. Let  $\theta^{*} \in \text{argmin}_{\theta \in \mathbb{R}^n} \overline{\mathcal{L}}(\theta)$ be a minimizer of the population risk $\overline{\mathcal{L}}(\theta)= \mathbb{E}_{Z_1^N}[\mathcal{L}(\theta;Z_1^N)]$. 

The task here is to estimate the unknown parameter $\theta^*$ from the data $Z_1^N$. A well known approach to this problem is to use a regularization technique, i.e., to solve the following convex optimization problem,
\begin{equation} \label{eq_c8Mregularization}
\hat{\theta}_{\lambda_N} \in \underset{\theta \in \mathbb{R}^n}{\text{arg min}} \ \ \lbrace \mathcal{L}(\theta; Z_1^N) + \lambda_N \mathcal{R}(\theta) \rbrace,
\end{equation}
where $\lambda_N > 0$ is a user-defined regularization parameter and $\mathcal{R}: \mathbb{R}^n \rightarrow \mathbb{R}^{+}$ is a norm. 

A difficulty when estimating the parameter $\theta^*$ using the above regularization technique is that one needs to find the regularization parameter $\lambda_N$. The traditional method to choose $\lambda_N$ is to use cross validation, i.e., to estimate the parameter $\theta^*$ with different values of $\lambda_N$, then select the value of $\lambda_N$ that provides the best fit to the validation data. This cross validation method is quite time consuming and very dependent on the data. Here we are specifically interested in the SPARSEVA (SPARSe Estimation based on a VAlidation criterion) framework, suggested in \cite{Rojas:11} and \cite{Rojas:14}, which provides automatic tuning of the regularization parameters.
Utilizing the SPARSEVA framework, an estimate of $\theta^*$ can be computed using the following convex optimization problem:
\begin{equation} \label{eq_c8sparseva}
\begin{aligned}
& \hat{\theta}_{\epsilon_N} \in \underset{\theta \in \mathbb{R}^n}{\text{arg min}}
& & \mathcal{R}(\theta) \\
& \ \ \ \ \ \ \ \ \ \text{s.t.}
& & \mathcal{L}(\theta; Z_1^N) \leq \mathcal{L}(\hat{\theta}_{NR}; Z_1^N)(1+\epsilon_N),
\end{aligned}
\end{equation}
where $\epsilon_N > 0$ is the regularization parameter and $\hat{\theta}_{NR}$ is the ``non-regularized'' estimate obtained from minimizing the cost function $\mathcal{L}(\theta; Z_1^N)$, i.e.
\begin{equation}
\hat{\theta}_{NR} \in \underset{\theta \in \mathbb{R}^n}{\text{arg min}} \ \ \ \ \mathcal{L}(\theta; Z_1^N).
\end{equation}
\CR{It can be shown \cite{Rojas:14} that \eqref{eq_c8Mregularization} and \eqref{eq_c8sparseva} are \emph{equivalent} in the sense that there exists a bijection between $\lambda_N$ and $\epsilon_N$ such that both estimators coincide. However, as discussed in \cite[Section V.D]{Rojas:14}, that bijection is data-dependent and it does not seem possible to derive an explicit expression for it.} The advantage of the SPARSEVA framework\CR{, with respect to \eqref{eq_c8Mregularization},} is that there are some natural choices of the regularization parameter $\epsilon_N$ based the chosen validation criterion. For example, as suggested in \cite{Rojas:11} \cite{Rojas:14}, $\epsilon_N$ can be chosen as $2n/N$ (Akaike Information Criterion (AIC)), $n\text{log}(N)/N$ (Bayesian Information Criterion (BIC)); or as suggested in \cite{Ha2015}, $n/N$ (Prediction Error Criterion).

For the traditional regularization method described in (\ref{eq_c8Mregularization}), \cite{Negahban2012} recently developed an upper bound on the estimation error between the estimate $\hat{\theta}_{\lambda_N}$ and the unknown parameter vector $\theta^*$. This bound is a function of some constants related to the nature of the data, the regularization parameter $\lambda_N$, the cost function $\mathcal{L}$ and the data length $N$. The beauty of this bound is that it quantifies the relationship between the estimation error and the finite data length $N$. Through this relationship, it is easy to confirm most of the properties of the estimate $\hat{\theta}_{\lambda_N}$ in the asymptotic scenario, i.e. $N \rightarrow \infty$, which were developed in the literature some time ago (\cite{Huang2008}, \cite{Knight2000}). 

Inspired by \cite{Negahban2012}, our goal is to derive a similar bound for the SPARSEVA estimate $\hat{\theta}_{\epsilon_N}$, i.e, we want to know how much the SPARSEVA estimate $\hat{\theta}_{\epsilon_N}$ differs from the true parameter $\theta^*$ when the data sample size $N$ is finite. Note that the notation and techniques used in this paper are similar to \cite{Negahban2012}; however, in \cite{Negahban2012}, the convex optimization problem is posed in the traditional regularization framework (\ref{eq_c8Mregularization}), while in this paper, the optimization problem is based on the SPARSEVA regularization (\ref{eq_c8sparseva}).

\section{Definitions and Properties of the Norm $\mathcal{R}(\theta)$ and the Cost Function $\mathcal{L}(\theta)$} 

In this section, we provide descriptions of some definitions and properties of the norm $\mathcal{R}(\theta)$ and the cost function $\mathcal{L}(\theta; Z_1^N)$, needed to establish an upper bound on the estimation error. Note that we only provide a brief summary such that the research described in this paper can be understood. Readers can find a more detailed discussion in \cite{Negahban2012}. 

\subsection{Decomposability of a Norm}

Let us consider a pair of \CR{arbitrary} linear subspaces \CR{of $\mathbb{R}^n$, $(\mathcal{M}, \overline{\mathcal{M}})$, such that $\mathcal{M} \subseteq \overline{\mathcal{M}}$}. The orthogonal complement of the space $\overline{\mathcal{M}}$ is then defined as,
\begin{equation} \nonumber
\overline{\mathcal{M}}^{\perp} = \lbrace v \in \mathbb{R}^n | \langle u,v\rangle = 0 \ \text{for all} \ u \in \overline{\mathcal{M}}\rbrace,
\end{equation}
where $\langle \cdot,\cdot \rangle$ is the inner product that maps $\mathbb{R}^n \times \mathbb{R}^n \rightarrow \mathbb{R}$. 

The norm $\mathcal{R}$ is said to be \textit{decomposable} with respect to $(\mathcal{M},\overline{\mathcal{M}}^{\perp})$ if
\begin{equation}
\mathcal{R}(\theta+\gamma) = \mathcal{R}(\theta) + \mathcal{R}(\gamma)
\end{equation}
for all $\theta \in \mathcal{M}$ and $\gamma \in \overline{\mathcal{M}}^{\perp}$.

There are many combinations of norms and vector spaces that satisfy this property (cf. \cite{Negahban2012}). An example is the $l_1$ norm and the sparse vector space defined (\ref{eq_normde}). For any subset $S \subseteq \lbrace 1,2,\dots,n\rbrace$ with cardinality $s$, define the model subspace $\mathcal{M}$ as,
\begin{equation} \label{eq_normde}
\mathcal{M}(S) = \lbrace\theta \in \mathbb{R}^n |\ \theta_j = 0\ \text{for all}\ j \not\in S\rbrace.
\end{equation}
Now if we define $\overline{\mathcal{M}}(S)=\mathcal{M}(S)$, then the orthogonal complement $\overline{\mathcal{M}}(S)$, with respect to the Euclidean inner product, can be computed as follows,
$$
\overline{\mathcal{M}}^{\perp}(S) = \lbrace\gamma \in \mathbb{R}^n | \gamma_j = 0\ \text{for all}\ j \in S\rbrace.
$$
As shown in \cite{Negahban2012}, the $l_1$-norm is decomposable with respect to the pair $(\mathcal{M}(S),\mathcal{M}^{\perp}(S))$.

\subsection{Dual Norm}

For a given inner product $\langle \cdot,\cdot \rangle$, the dual of the norm $\mathcal{R}$ is defined by,
\begin{equation} \label{eq:dual_norm}
\mathcal{R}^*(v) = \underset{u \in \mathbb{R}^n \setminus \lbrace 0 \rbrace}{\text{sup}} \dfrac{\langle u,v \rangle}{\mathcal{R}(u)} = \underset{\mathcal{R}(u) \leq 1}{\text{sup}} \langle u, v \rangle,
\end{equation}
where $\text{sup}$ is the supremum operator. \\
Based on the above definition, one can easily see that the dual of the $l_1$ norm, with respect to the Euclidean inner product, is the $l_{\infty}$ norm \cite{Negahban2012}.

\subsection{Strong Convexity}

A twice differentiable function $\mathcal{L}(\theta): \mathbb{R}^n\rightarrow\mathbb{R}$ is \emph{strongly convex} on $\mathbb{R}^n$ when there exists an $m>0$ such that its Hessian $\bigtriangledown^2\mathcal{L}(\theta)$ satisfies,
\begin{equation} \label{eq_c8SC}
\bigtriangledown^2\mathcal{L}(\theta) \succeq mI,
\end{equation}
for all $\theta \in \mathbb{R}^n$ \cite{Boyd2004}. %When $\mathcal{L}(\theta)$ is twice continuously differentiable, then 
This is equivalent to the statement that the minimum eigenvalue of $\bigtriangledown^2\mathcal{L}(\theta)$ is not smaller than $m$ for all $\theta \in \mathbb{R}^n$. \\
An interesting consequence of the strong convexity property in (\ref{eq_c8SC}) is that for all $\theta, \Delta \in \mathbb{R}^n$, we have,
\begin{equation} \label{eq_kl_def}
\mathcal{L}(\theta+\Delta) \geq \mathcal{L}(\theta) + \bigtriangledown\mathcal{L}(\theta)^T\Delta+\dfrac{m}{2}\Vert \Delta \Vert_2^2.
\end{equation}
The inequality in (\ref{eq_kl_def}) has a geometric interpretation in that the graph of the function $\mathcal{L}(\theta)$ has a positive curvature at any $\theta \in \mathbb{R}^n$. The term $m/2$ for the largest $m$ satisfying (\ref{eq_c8SC}) is typically known as the \emph{curvature} of $\mathcal{L}(\theta)$. 

\subsection{Subspace Compatibility Constant}

For a given norm $\mathcal{R}$ and an error norm $\Vert \cdot \Vert$, the \emph{subspace compatibility constant} of a subspace \CR{$\mathcal{M} \subseteq \mathbb{R}^n$} with respect to the pair $(\mathcal{R}, \Vert \cdot \Vert)$ is defined as,
\begin{equation} \label{Sub_CompatC}
\Psi(\mathcal{M}) = \underset{u \in \mathcal{M}\setminus\lbrace0\rbrace}{\sup} \dfrac{\mathcal{R}(u)}{\Vert u \Vert}.
\end{equation}
This quantity measures how well the norm $\mathcal{R}$ is compatible with the error norm $\Vert . \Vert$ over the subspace $\mathcal{M}$. As shown in \cite{Negahban2012}, when $\mathcal{M}$ is $\mathbb{R}^s$, the regularized norm $\mathcal{R}$ is the $l_1$ norm, and the error norm is the $l_2$ norm, then the subspace compatibility constant is $\Psi(\mathcal{M}) = \sqrt{s}$.
Notice also that $\Psi(\mathcal{M})$ is finite, due to the equivalence of finite dimensional norms.

\subsection{Projection Operator} \label{subsec:projection}

The projection of a vector $u$ onto a space $\mathcal{M}$, with respect to the Euclidean norm, is defined by the following,
\begin{equation}
\Pi_{\mathcal{M}}(u) = \argmin\limits_{v\in\mathcal{M}} \Vert u-v\Vert_2.
\end{equation}
In the sequel, to simplify the notation, we will write $u_{\mathcal{M}}$ to denote $\Pi_{\mathcal{M}}(u)$.

\section{Analysis of the Regularization Technique using the SPARSEVA}

In this section, we apply the properties described in Section 3 to derive an upper bound on the error between the SPARSEVA estimate $\hat{\theta}_{\epsilon_N}$ and the unknown parameter $\theta^*$. This upper bound is described in the following theorem.

\begin{theorem} \label{tr_generalbound}
Assume $\mathcal{R}$ is a norm and is decomposable with respect to the subspace pair ($\mathcal{M},\overline{\mathcal{M}}^\perp$) and the cost function $\mathcal{L}(\theta)$ is differentiable and strongly convex with curvature $\kappa_{\mathcal{L}}$. Consider the SPARSEVA problem in (\ref{eq_c8sparseva}), then the following properties hold:
\begin{enumerate} [i.]
\item When $\epsilon_N > 0$, there exists a Lagrange multiplier, $\lambda_N = \lambda_{\epsilon_N}$, such that \eqref{eq_c8Mregularization} and \eqref{eq_c8sparseva} have the same solution.
%for the SPARSEVA problem in (\ref{eq_c8sparseva}).
\item Any optimal solution $\hat{\theta}_{\epsilon_N} \neq 0$ of the SPARSEVA problem (\ref{eq_c8sparseva}) satisfies the following inequalities: 
\vspace{0.1cm} \\
\begin{itemize}
\item If $\epsilon_N$ is chosen such that $$\lambda_{\epsilon_N} \leq 1/\mathcal{R}^*(\nabla\mathcal{L}(\theta^*)),$$ then
\begin{equation} \label{bound_gn_case1}
\begin{aligned}
\Vert \hat{\theta}_{\epsilon_N}-\theta^* \Vert^2_2\ & \leq \ \dfrac{4}{\kappa^2_{\mathcal{L}}\lambda^2_{\epsilon_N}}\Psi^2(\overline{\mathcal{M}}) + \dfrac{4}{\kappa_{\mathcal{L}}\lambda_{\epsilon_N}}\mathcal{R}(\theta^*_{\mathcal{M}^{\perp}}).
\end{aligned}
\end{equation}
\item If $\epsilon_N$ is chosen such that $$\lambda_{\epsilon_N} > 1/\mathcal{R}^*(\nabla\mathcal{L}(\theta^*)),$$ then
\begin{equation} \label{bound_gn_case2}
\begin{aligned}
\Vert \hat{\theta}_{\epsilon_N}-\theta^* \Vert^2_2\ & \leq \ \dfrac{2}{\kappa^2_{\mathcal{L}}}\lbrace\mathcal{R}^*(\nabla\mathcal{L}(\theta^*))\rbrace^2 \Psi^2(\overline{\mathcal{M}}) \\
& \ \ \ +\dfrac{8}{\kappa^2_{\mathcal{L}}}\lbrace\mathcal{R}^*(\nabla\mathcal{L}(\theta^*))\rbrace^2 \Psi^2(\overline{\mathcal{M}}^\perp) \\
& \ \ \ + \dfrac{4}{\kappa_{\mathcal{L}}\lambda_{\epsilon_N}}\mathcal{R}(\theta^*_{\mathcal{M}^{\perp}}).
\end{aligned}
\end{equation}
\end{itemize}
\end{enumerate}
\end{theorem}

\textit{\textbf{Proof.}} See the Appendix (Section A.2). \hfill $\square$

\begin{remark}
Note that Theorem \ref{tr_generalbound} is intended to provide an upper bound on the estimation error for the general SPARSEVA problem (\ref{eq_c8sparseva}). At this stage, it is hard to evaluate, or quantify, the value on the right hand side of the inequalities (\ref{bound_gn_case1}) and (\ref{bound_gn_case2}) as they still contain the term $\lambda_{\epsilon_N}$ and other abstract terms. However, in the later sections of this paper, from this general upper bound, we will provide bounds on the estimation errors for some specific scenarios.
\end{remark}

\begin{remark}
The bound in Theorem \ref{tr_generalbound} is actually a family of bounds. For each choice of the pair of subspaces ($\mathcal{M},\overline{\mathcal{M}}^\perp$), there is one bound for the estimation error. Hence, in the usual sense, to apply Theorem \ref{tr_generalbound} for any specific scenario, the goal is to choose $\mathcal{M}$ and $\overline{\mathcal{M}}^\perp$ to obtain an optimal rate of the bound.
\end{remark}

%%%%%%%%%%%%%%%%%%%%%%%%%%%%%%%%%%%%%%%%%%%%%%
\section{An Upper Bound for Sparse Regression}
%%%%%%%%%%%%%%%%%%%%%%%%%%%%%%%%%%%%%%%%%%%%%%

In this section, we illustrate how to apply Theorem \ref{tr_generalbound} to derive an upper bound of the error between the SPARSEVA estimate $\hat{\theta}_{\epsilon_N}$ and the true parameter $\theta^*$ for the following linear regression model,
\begin{equation} \label{eq_lr}
Y_N = \Phi_N^T\theta^* + e,
\end{equation}
where $\theta^* \in \mathbb{R}^n$ is the unknown parameter that is required to be estimated; $e \in \mathbb{R}^N$ is the disturbance noise; $\Phi_N \in \mathbb{R}^{n\times N}$ is the regression matrix and $Y_N\in \mathbb{R}^{N}$ is the output vector. Here, we make the following assumption on the true parameter $\theta^*$,

\begin{assumption} \label{assu_theta}
The true parameter $\theta^*$ is ``weakly" sparse, i.e. $\theta^* \in \mathbb{B}_q(R_q)$, where,
\begin{equation} \label{eq_BallRq}
\mathbb{B}_q(R_q) := \left\lbrace \theta \in \mathbb{R}^p \left\vert \sum_{i=1}^p |\theta_i|^q \leq R_q \right.\right\rbrace,
\end{equation}
with $q \in [0,1]$ being a constant. 
\end{assumption}

Using the SPARSEVA framework in (\ref{eq_c8sparseva}) with $\mathcal{R}$ chosen as the $l_1$ norm and the cost function $\mathcal{L}(\theta)$ chosen as,
\begin{equation} \label{eq_costFunc}
\mathcal{L}(\theta) = \dfrac{1}{2N}\Vert Y_N - \Phi_N^T\theta \Vert_2^2,
\end{equation}
then an estimate of $\theta^*$ in (\ref{eq_lr}) can be found by solving the following problem,
\begin{equation} \label{sparseva_ln}
\begin{aligned}
& \hat{\theta}_{\epsilon_N} \in \underset{\theta \in \mathbb{R}^n}{\text{arg min}}
& & \Vert\theta\Vert_1 \\
& \ \ \ \ \ \ \ \ \ \ \text{s.t.}
& & \mathcal{L}(\theta) \ \leq\ \mathcal{L}(\hat{\theta}_{NR})(1+\epsilon_N),
\end{aligned}
\end{equation}
with $\hat{\theta}_{NR}=(\Phi_N\Phi_N^T)^{-1}\Phi_NY_N$ and $\epsilon_N > 0$ being the user-defined regularization parameter. Now $\epsilon_N$ can be chosen as either $2n/N$ or $\text{log}(N)n/N$ as suggested in \cite{Rojas:11}; or $n/N$ as suggested in \cite{Ha2015}.

\begin{remark}
Note that the sparse regression problem is very common in system identification and is often used to obtain a low order linear model by regularization.
\end{remark}

\begin{remark}
For Assumption \ref{assu_theta}, note that when $q=0$, under the convention that $0^0 = 0$, the set in (\ref{eq_BallRq}) corresponds to an exact sparsity set, where all the elements belonging to the set have at most $R_0$ non-zero entries. Generally, for $q \in (0,1]$, the set $\mathbb{B}_q(R_q)$ forces the ordered absolute values of $\theta^*$ to decay with a certain rate.
\end{remark}

\subsection{An Analysis on the Strong Convexity Property and the Curvature of the $l_2$ norm Cost Function}

Consider the convex optimization problem in (\ref{sparseva_ln}), the Hessian matrix of the cost function $\mathcal{L}(\theta)$ is computed as,
$$\bigtriangledown^2\mathcal{L}(\theta) = \dfrac{1}{N}\Phi_N\Phi_N^T.$$
To prove that $\mathcal{L}(\theta)$ is strongly convex, we need to prove,
\begin{equation} \label{eq_sc}
\exists \kappa_{\mathcal{L}} >0 \ \ \ \text{s.t.} \ \ \ \dfrac{1}{N} \Phi_N\Phi_N^T \succeq 2\kappa_{\mathcal{L}}I.
\end{equation}
We see that the requirement in (\ref{eq_sc}) coincides with the requirement of persistent excitation of the input signal in a system identification problem. If an experiment is well-designed, then the input signal $u(t)$ needs to be persistently exciting of order $n$, i.e., the matrix $\Phi_N\Phi_N^T$ is a positive definite matrix. This means that the condition in (\ref{eq_sc}) is always satisfied for any linear regression problem derived from a well posed system identification problem. This means that for any choice of the regression matrix $\Phi_N$ that satisfies the persistent excitation condition, there exists a positive curvature $\kappa_{\mathcal{L}}$ of the cost function $\mathcal{L}(\theta)$.

Consider $\Phi_N \in \mathbb{R}^{n\times N}$ to be a matrix where each row $\Phi_{N,j}$ is sampled from a Normal distribution of zero mean and covariance matrix $\Sigma \in \mathbb{R}^{N\times N}$, i.e., $\Phi_{N,j} \sim \mathcal{N}(0,\Sigma), \ \ \ \forall j=1,..,n.$ We then denote the distribution of the smallest eigenvalue of $N^{-1} \Phi_N\Phi_N^T$ to be $P(x\vert\Sigma,N,n)$, means that given a probability $1-\alpha\ ,\ 0\leq \alpha \leq 1$, there exists a value $w_{\text{min}}$ such that $N^{-1} \Phi_N\Phi_N^T \succeq w_{\text{min}}I$, for any matrix $\Phi_N$ constructed following the above assumption. Then the global curvature $\kappa$, i.e. the curvature that satisfies (\ref{eq_sc}) for any regression matrix $\Phi_N$, can be expressed as $(1/2) w_{\text{min}}$. For the rest of the paper, we will denote by $\kappa_\alpha$ lower bound on the global curvature $\kappa$ with probability $1-\alpha \ ,\ 0\leq \alpha \leq 1$.

\subsection{Assumptions}

For the linear regression in (\ref{eq_lr}), the following assumptions are made:

\begin{assumption} \label{assu_X}
The rows $\Phi_{N,j}, j=1,...n$ of the regressor matrix $\Phi_N$ are distributed as $\Phi_{N,j} \sim \mathcal{N}(0, \Sigma)$, where $\Sigma \in \mathbb{R}^{N \times N}$ is a constant, symmetric, positive definite matrix.
\end{assumption}
Note that an obvious practical case where Assumption \ref{assu_X} is satisfied is when the model is FIR and the input signal being white noise or coloured noise.

\begin{assumption} \label{assu_e}
The noise vector $\text{e} \in \mathbb{R}^N$ is Gaussian with i.i.d. $\mathcal{N}(0,\sigma_e^2)$ entries.\footnote{\CR{The assumption of Gaussian noise is fairly standard in system identification. However, this assumption can be relaxed to `sub-Gaussian' noise (\emph{i.e.}, when the tails of the noise distribution decay like $e^{-\alpha x^2}$) at the expense of longer derivations.}}
\end{assumption}

\subsection{Developing the Upper Bound}

\CR{The following theorem provides an upper bound on the estimation error $\Vert \hat{\theta}_{\epsilon_N} - \theta^* \Vert_2$ for the optimization problem in (\ref{sparseva_ln}) in the case of weakly sparse estimates.}

\begin{theorem} \label{tr_spbound}
Suppose Assumptions \ref{assu_X}, \ref{assu_e} and \ref{assu_theta} hold, when N is large, then with probability $(1-\alpha)(1 - 4 n \beta)$ ($0 \leq \alpha \leq 1,\ 0 \leq \beta \leq 1$), if $\hat{\theta}_{\epsilon_N} \neq 0$ we have the following inequality \\
\begin{equation} \label{eq_tr_spbound}
\Vert \hat{\theta}_{\epsilon_N} - \theta^* \Vert^2_2\ \leq \max(a_1, a_2),
\end{equation}
where
\begin{align*}
a_1 & = \frac{8 n_\eta \sigma_e^2 s_{\text{max}} \chi^2_{\beta}(N-n) (1+\epsilon_N) \ln(2/\beta)}{\kappa_{\alpha}^2 N^2} \\
&\quad + \frac{\sqrt{32\sigma_e^2 s_{\text{max}} \chi^2_{\beta}(N-n) (1+\epsilon_N)\ln(2/\beta)}}{\kappa_{\alpha} N} \Vert \theta_{[n_{\eta}+1:n]}^* \Vert_1, \\
a_2 & = \frac{(16 n - 12 n_\eta) \sigma_e^2 \chi^2_{\beta}(\Sigma, I) \ln(2/\beta)}{\kappa_{\alpha}^2 N^2} \\
& \ \ \ \ +\frac{\sqrt{32\sigma_e^2 s_{\text{max}} \chi^2_{\beta}(N-n) (1+\epsilon_N)\ln(2/\beta)}}{\kappa_{\alpha} N} \Vert \theta_{[n_{\eta}+1:n]}^* \Vert_1.
\end{align*}
where $\kappa_\alpha$ is a lower bound on the curvature of the regression matrix (i.e., half the smallest eigenvalue of $N^{-1}\Phi_N^T \Phi_N$) with probability $1-\alpha$, $n_{\eta}$ is any integer between $1$ and $n$, $\theta_{[n_{\eta}+1:n]}^*$ is the vector formed from the $n-n_{\eta}$ smallest (in magnitude) entries of $\theta^*$, and $s_{\text{max}}$ is the maximum singular value of the matrix $\Sigma$.
\end{theorem}

\vspace{0.1cm}
\textit{\textbf{Proof.}} \CR{This proof relies on three preliminary results introduced in Appendix A.3.} For an integer $n_{\eta} \in \{1,\dots,n\}$, define $S_{\eta}$ as the set of the indices of the $n_{\eta}$ largest (in magnitude) entries of $\theta^*$, and its complementary set $S^c_{\eta}$ as
\begin{equation}
\begin{aligned}
& S^c_{\eta} = \lbrace 1,2,...,n \rbrace \setminus S_{\eta};
\end{aligned}
\end{equation}
with the corresponding subspaces $\mathcal{M}(S_\eta)$ and $\mathcal{M}^\perp(S_\eta)$ as,
\begin{equation}
\begin{aligned}
\mathcal{M}(S_\eta) &= \lbrace \theta \in \mathbb{R}^n \mid \theta_j=0 \ \ \forall j \not\in S_\eta \rbrace, \\
\mathcal{M}^\perp(S_\eta) &= \lbrace \gamma \in \mathbb{R}^n \mid \gamma_j=0 \ \ \forall j \in S_\eta \rbrace.
\end{aligned}
\end{equation}
Using the definition of the subspace compatibility constant described in Section 3, we have, 
\begin{equation} \label{eq_Scardian}
\begin{aligned}
\Psi^2(\mathcal{M}(S_\eta))=\vert S_\eta \vert = n_{\eta}, \\
\Psi^2(\mathcal{M}^\perp(S_\eta)) = \vert S^c_\eta \vert = n - n_{\eta}.
\end{aligned}
\end{equation}
where $\vert S \vert$ denotes the cardinality of $S$.

Now, for Theorem \ref{tr_generalbound} to generate an upper bound for the problem (\ref{sparseva_ln}), we need to establish an upper bound on $\Vert \theta^*_{\mathcal{M}^{\perp}(S_{\eta})} \Vert_1$. Based on the definition of the subspace $\mathcal{M}^\perp(S_\eta)$, we have,
\begin{equation} \label{eq_Rthetastar}
\Vert \theta^*_{\mathcal{M}^{\perp}(S_{\eta})} \Vert_1 = \Vert \theta_{[n_{\eta}+1:n]}^* \Vert_1,
\end{equation}
where $\theta_{[n_{\eta}+1:n]}^*$ denotes the vector formed from the $n-n_{\eta}$ smallest (in magnitude) entries of $\theta^*$. Define $\kappa_{\alpha}$ as a lower bound on the global curvature of the regression matrix, i.e. half the smallest eigenvalue of $\Phi_N^T\Phi_N$, with probability $1-\alpha,\ 0 \leq \alpha \leq 1$. Substituting \CR{the results of Propositions~\ref{pr_Lagrange}-\ref{pr_derThetahat} from Appendix A.3}, (\ref{eq_Scardian}) and (\ref{eq_Rthetastar}) %(\ref{eq_lagrange}), (\ref{eq_pr3spc}), (\ref{eq_pr4spc})
into the bound in Theorem \ref{tr_generalbound}, then with $n_{\eta}$ being any integer between $1$ and $n$, we have the following bounds:

\begin{itemize}
\item
If $\epsilon_N$ is chosen such that
$$\lambda_{\epsilon_N} \leq 1/\mathcal{R}^*(\nabla\mathcal{L}(\theta^*)) \ = 1/ \| \nabla\mathcal{L}(\theta^*)\|_\infty,$$
then, with probability at least $(1-\alpha)(1 - 2 n \beta)$,
\begin{align*}
& \Vert \hat{\theta}_{\epsilon_N}-\theta^* \Vert^2_2\ \\%
&\leq \frac{8 n_\eta \sigma_e^2 s_{\text{max}} \chi^2_{\beta}(N-n) (1+\epsilon_N) \ln(2/\beta)}{\kappa_{\alpha}^2 N^2} \\
&\quad + \frac{\sqrt{32\sigma_e^2 s_{\text{max}} \chi^2_{\beta}(N-n) (1+\epsilon_N)\ln(2/\beta)}}{\kappa_{\alpha} N} \Vert \theta_{[n_{\eta}+1:n]}^* \Vert_1.
\end{align*}
\item
If $\epsilon_N$ is chosen such that
$$\lambda_{\epsilon_N} > 1/\mathcal{R}^*(\nabla\mathcal{L}(\theta^*)) \ = 1/ \| \nabla\mathcal{L}(\theta^*)\|_\infty,$$
then, with probability at least $(1-\alpha)(1 - 4 n \beta)$,
\begin{align*}
&\Vert \hat{\theta}_{\epsilon_N}-\theta^* \Vert^2_2 \\ %
& \leq \frac{(16 n - 12 n_\eta) \sigma_e^2 \chi^2_{\beta}(\Sigma, I) \ln(2/\beta)}{\kappa_{\alpha}^2 N^2} \\
& \ \ \ \ +\frac{\sqrt{32\sigma_e^2 s_{\text{max}} \chi^2_{\beta}(N-n) (1+\epsilon_N)\ln(2/\beta)}}{\kappa_{\alpha} N} \Vert \theta_{[n_{\eta}+1:n]}^* \Vert_1.
\end{align*}
\end{itemize}
\vspace{-0.1cm}
Therefore, for $n_{\eta}$ being any integer between $1$ and $n$, with probability at least ($1-\alpha$)($1 - 4n \beta)$, we have
\begin{equation}
\Vert \hat{\theta}_{\epsilon_N}-\theta^* \Vert^2_2\ \leq \text{max}(a_1, a_2),
\end{equation}
where
\begin{align*}
a_1 & = \frac{8 n_\eta \sigma_e^2 s_{\text{max}} \chi^2_{\beta}(N-n) (1+\epsilon_N) \ln(2/\beta)}{\kappa_{\alpha}^2 N^2} \\
&\quad + \frac{\sqrt{32\sigma_e^2 s_{\text{max}} \chi^2_{\beta}(N-n) (1+\epsilon_N)\ln(2/\beta)}}{\kappa_{\alpha} N} \Vert \theta_{[n_{\eta}+1:n]}^* \Vert_1,
\end{align*}
\begin{align*}
a_2 & = \frac{(16 n - 12 n_\eta) \sigma_e^2 \chi^2_{\beta}(\Sigma, I) \ln(2/\beta)}{\kappa_{\alpha}^2 N^2} \\
& \ \ \ \ +\frac{\sqrt{32\sigma_e^2 s_{\text{max}} \chi^2_{\beta}(N-n) (1+\epsilon_N)\ln(2/\beta)}}{\kappa_{\alpha} N} \Vert \theta_{[n_{\eta}+1:n]}^* \Vert_1.
\end{align*}
\hfill $\square$

\begin{remark}
The bound in Theorem \ref{tr_spbound} is also a family of bounds, one for each value of $n_{\eta}$. 
\end{remark}

\begin{remark}
When $\Sigma = \sigma_u I$, i.e. the model is FIR and the input $u(t)$ is white noise, then $s_{\text{max}}=\sigma_u$ and the generalized Chi square distribution $\chi^2(\Sigma, I)$ becomes the Chi square distribution $\sigma_u\chi^2(N)$.
\end{remark}

\begin{remark}
Note that the developed bound in Theorem \ref{tr_spbound} depends on the true parameter $\theta^*$, which is unknown but constant. Using a similar proof as in Proposition 2.3 of \cite{Foucart:04}, we can derive under Assumption \ref{assu_theta} an upper bound for the term $\Vert \theta_{[n_{\eta}+1:n]}^* \Vert_1$. Specifically, we have,
\begin{align*}
\Vert \theta_{[n_{\eta}+1:n]}^* \Vert_1 = \sum_{i=n_\eta+1}^n |\theta^*_{[i]}| = \sum_{i=n_\eta+1}^n |\theta^*_{[i]}|^{1-q} |\theta^*_{[i]}|^q
\end{align*}
Since $S_{\eta}$ is the set of the indices of the $n_{\eta}$ largest (in magnitude) entries of $\theta^*$, i.e. $|\theta^*_{[i]}| \leq |\theta^*_{n_\eta}|,\ \forall i = n_\eta+1, ...,n$, hence,
\begin{align*}
\Vert \theta_{[n_{\eta}+1:n]}^* \Vert_1 & \leq |\theta^*_{n_\eta}|^{1-q} \sum_{i=n_\eta+1}^n |\theta^*_{[i]}|^q
\end{align*}
Using the same argument, we have,
\begin{align*}
|\theta^*_{n_\eta}|^{1-q} = \Big(\dfrac{1}{n_\eta} \sum_{i=1}^{n_\eta} |\theta^*_{[n_\eta]}|^q \Big)^{(1-q)/q} \leq \Big(\dfrac{1}{n_\eta} \sum_{i=1}^{n_\eta} |\theta^*_{[i]}|^q \Big)^{(1-q)/q}.
\end{align*}
Therefore,
\begin{align*}
\Vert \theta_{[n_{\eta}+1:n]}^* \Vert_1 & \leq \Big(\dfrac{1}{n_\eta} \sum_{i=1}^{n_\eta} |\theta^*_{[i]}|^q \Big)^{(1-q)/q} \sum_{i=n_\eta+1}^n |\theta^*_{[i]}|^q \\
& \leq \Big(\dfrac{1}{n_\eta} \sum_{i=1}^{n} |\theta^*_{[i]}|^q \Big)^{(1-q)/q} \sum_{i=1}^n |\theta^*_{[i]}|^q \\
& \leq \Big(\dfrac{1}{n_\eta}\Big)^{1/q-1} \Vert \theta^* \Vert_q^{1-q} \Vert \theta^* \Vert_q^{q}\\
& \leq (n_\eta)^{1-1/q} \Vert \theta^* \Vert_q \\
& \leq (n_\eta)^{1-1/q} (R_q)^{1/q}.
\end{align*}
This means we can always place an upper bound on the term $\Vert \theta_{[n_{\eta}+1:n]}^* \Vert_1$ by a known constant which depends on the nature of the true parameter $\theta^*$. Therefore, from Theorem \ref{tr_spbound}, we can see that the estimation error $\Vert \hat{\theta}_{\epsilon_N}-\theta^* \Vert^2_2 = O_p(N^{-1/2})$ \cite{Rojas:14}. This confirms the result in \cite{Rojas:14}, that in the asymptotic case, when $\epsilon_N > 0$, the SPARSEVA estimate $\hat{\theta}_{\epsilon_N}$ converges to the true parameter $\theta^*$. 
%
%This also means that the value of $\Vert \theta_{[n_{\eta}+1:n]}^* \Vert_1$ is small, which guarantees the tightness of the developed bound.
\end{remark}

%%%%%%%%%%%%%%%%%%%%%%%
\section{Numerical Evaluation}
%%%%%%%%%%%%%%%%%%%%%%%

In this section, numerical examples are presented to illustrate the bound $ \Vert \hat{\theta}_{\epsilon_N}-\theta^* \Vert^2_2$ as stated in Theorem \ref{tr_spbound}. In Section 6.1, we consider the case when the input is Gaussian white noise whilst in Section 6.2, the input is a correlated signal with zero mean.

\subsection{Gaussian White Noise Input}

In this section, a random discrete time system with a random model order between 1 and 10 is generated using the command \textit{drss} from Matlab. The system has poles with magnitude less than 0.9. Gaussian white noise is added to the system output to give different levels of SNR, e.g. 30dB, 20dB and 10dB. For each noise level, 50 different input excitation signals (Gaussian white noise with variance 1) and output noise realizations are generated. For each set of input and output data, the system parameters are estimated using a different sample size, i.e., $N=[450, 1000, 5000, 10000, 50000, 100000]$. 

The FIR model structure is used here in order to construct the SPARSEVA problem (\ref{sparseva_ln}). The number of parameters $n$ of the FIR model is set to be 35. The regularization parameter $\epsilon_N$ is chosen as $n/N$ \cite{Ha2015}.

We then compute the upper bound of $\Vert \hat{\theta}_{\epsilon_N}-\theta^* \Vert_2$ using (\ref{eq_tr_spbound}) with different values of $n_{\eta}$, i.e. $n_{\eta} = [10, 15, 25]$. The probability parameters $\alpha$ and $\beta$ are chosen to be $0.02$ and $0.001$ respectively. Related to the computation of the universal constant $\kappa$ corresponding to the distribution $\mathcal{N}(0,\Sigma)$, note that, in reality, it is very difficult to compute its exact distribution $P(x\vert\Sigma,N,n)$, hence, here we use an empirical method to compute the distribution $P(x\vert\Sigma,N,n)$. The idea is to generate a large number of random matrices $\Phi_N$, compute the smallest eigenvalue of $N^{-1} \Phi_N\Phi_N^T$, and then build a histogram of these values, which is an approximation of $P(x\vert\Sigma,N,n)$. Then we compute the value of $w_{\text{min}}$ to ensure the inequality $N^{-1} \Phi_N\Phi_N^T \succeq w_{\text{min}}I$ occurs with probability $1-\alpha$. Finally, $\kappa_\alpha$ is computed using the formula $\kappa_\alpha = (1/2) w_{\text{min}}$.

With the setting described above, the probability of the upper bound being correct is $(1-\alpha)(1-4n\beta)=0.84$. This upper bound will be compared with $\Vert \hat{\theta}_{\epsilon_N}-\theta^* \Vert_2$. Note that we plot both the upper bound and the true estimation errors on a logarithmic scale.

\begin{figure}[h]
\begin{center}
\vspace{0.2cm}
\includegraphics[width=79mm,height=53mm]{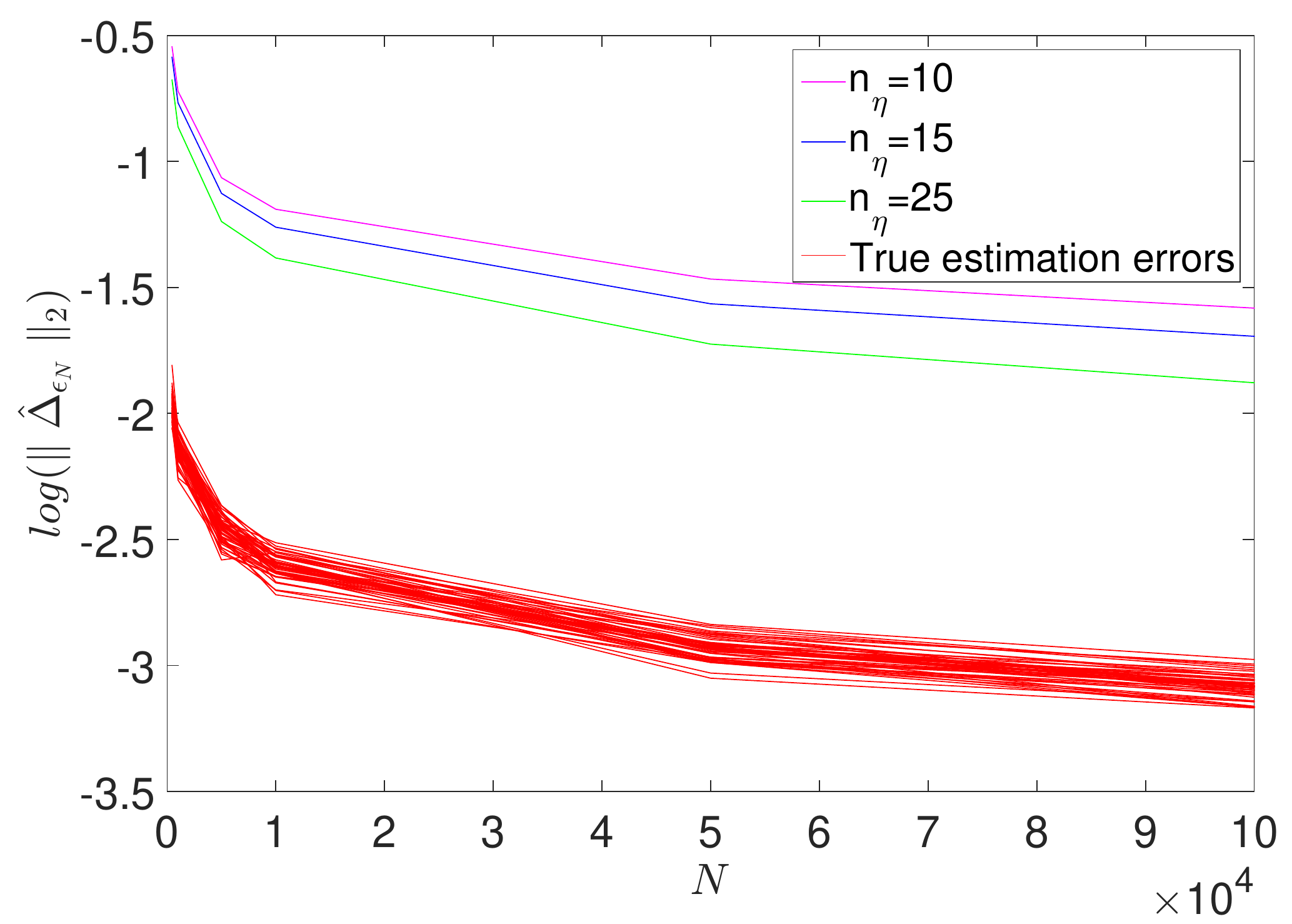}
\caption{Plot of the estimation error for a SNR=30dB and a Gaussian white input signal.} 
\label{fig_wn30dB} 
\end{center}
\end{figure}

\begin{figure}[h]
\begin{center}
\vspace{0.2cm}
\includegraphics[width=77mm,height=53mm]{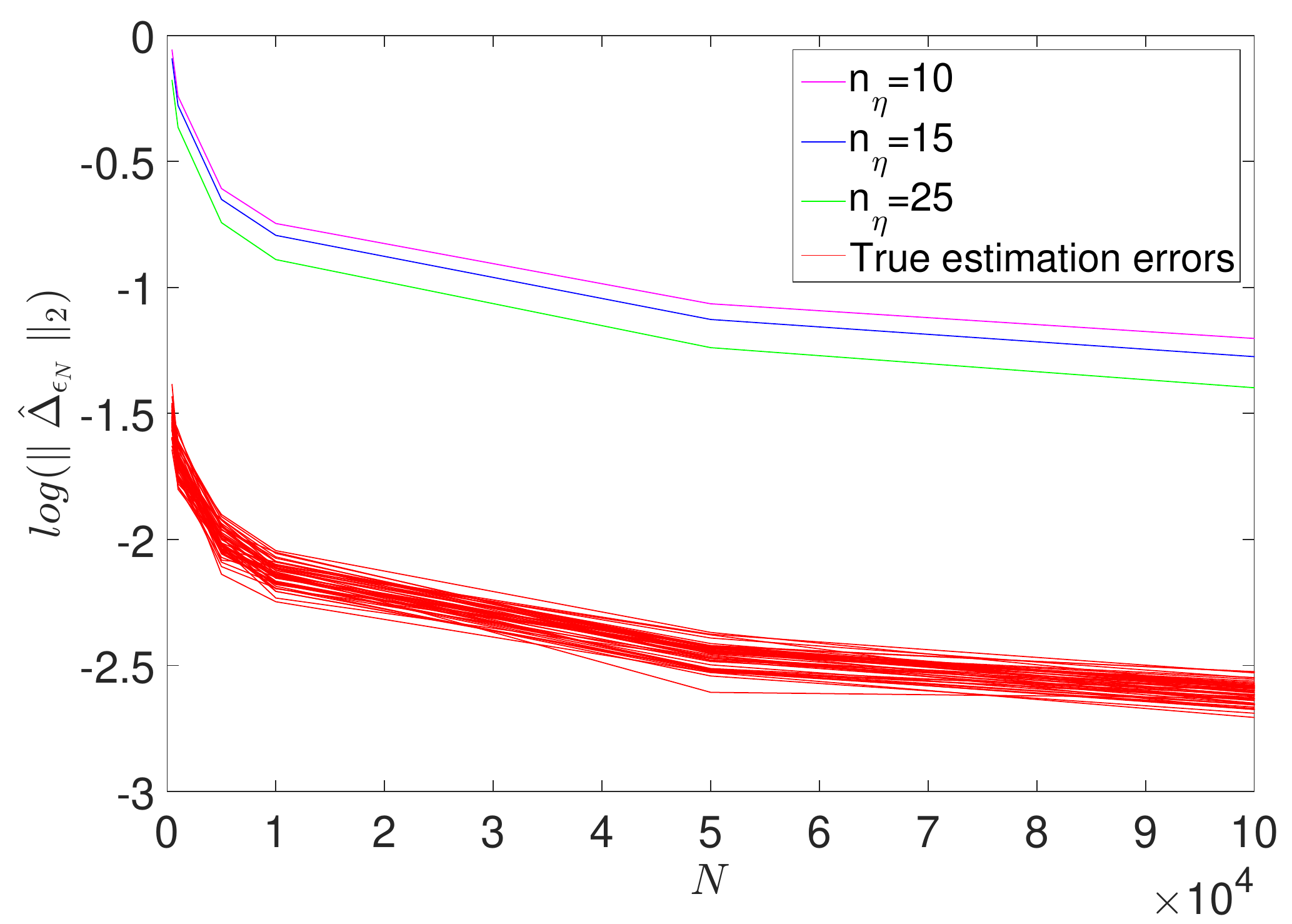}
\caption{Plot of the estimation error for a SNR=20dB and a Gaussian white input signal.} 
\label{fig_wn20dB} 
\end{center}
\end{figure}

\begin{figure}[h]
\begin{center}
\vspace{0.2cm}
\includegraphics[width=77mm,height=53mm]{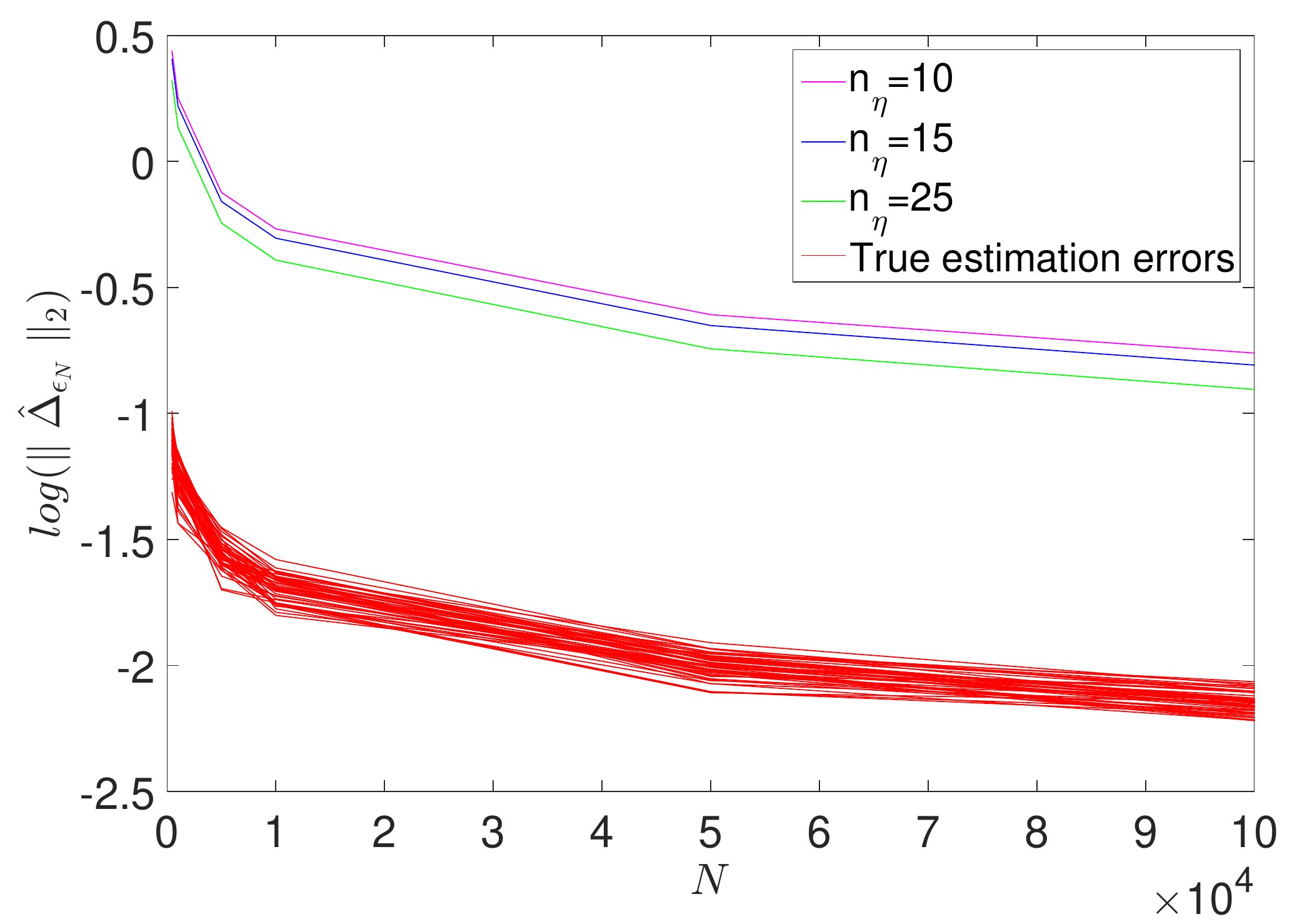}
\caption{Plot of the estimation error for a SNR=10dB and a Gaussian white input signal.} 
\label{fig_wn10dB} 
\end{center}
\end{figure}

Plots of the estimation error versus the data length $N$ with different noise levels are displayed in Figures \ref{fig_wn30dB} to \ref{fig_wn10dB}. In Figures \ref{fig_wn30dB} to \ref{fig_wn10dB}, the red lines are the true estimation errors from 50 estimates using the SPARSEVA framework. The magenta, blue and cyan lines are the upper bounds developed in Theorem \ref{tr_spbound}, which correspond to $n_\eta = [10,15,25]$, respectively. We can see that the plots confirm the bound developed in Theorem \ref{tr_spbound} for all noise levels. When $N$ becomes large, the estimation error and the corresponding upper bound become smaller. When $N$ goes to infinity, the estimation error will tend to 0. Note that the bounds are slightly different for the chosen values of $n_\eta$, however, not significantly. As can be seen, the bounds are relatively insensitive to the choice of $n_\eta$.

\begin{figure}[h]
\begin{center}
\vspace{0.2cm}
\includegraphics[width=77mm,height=53mm]{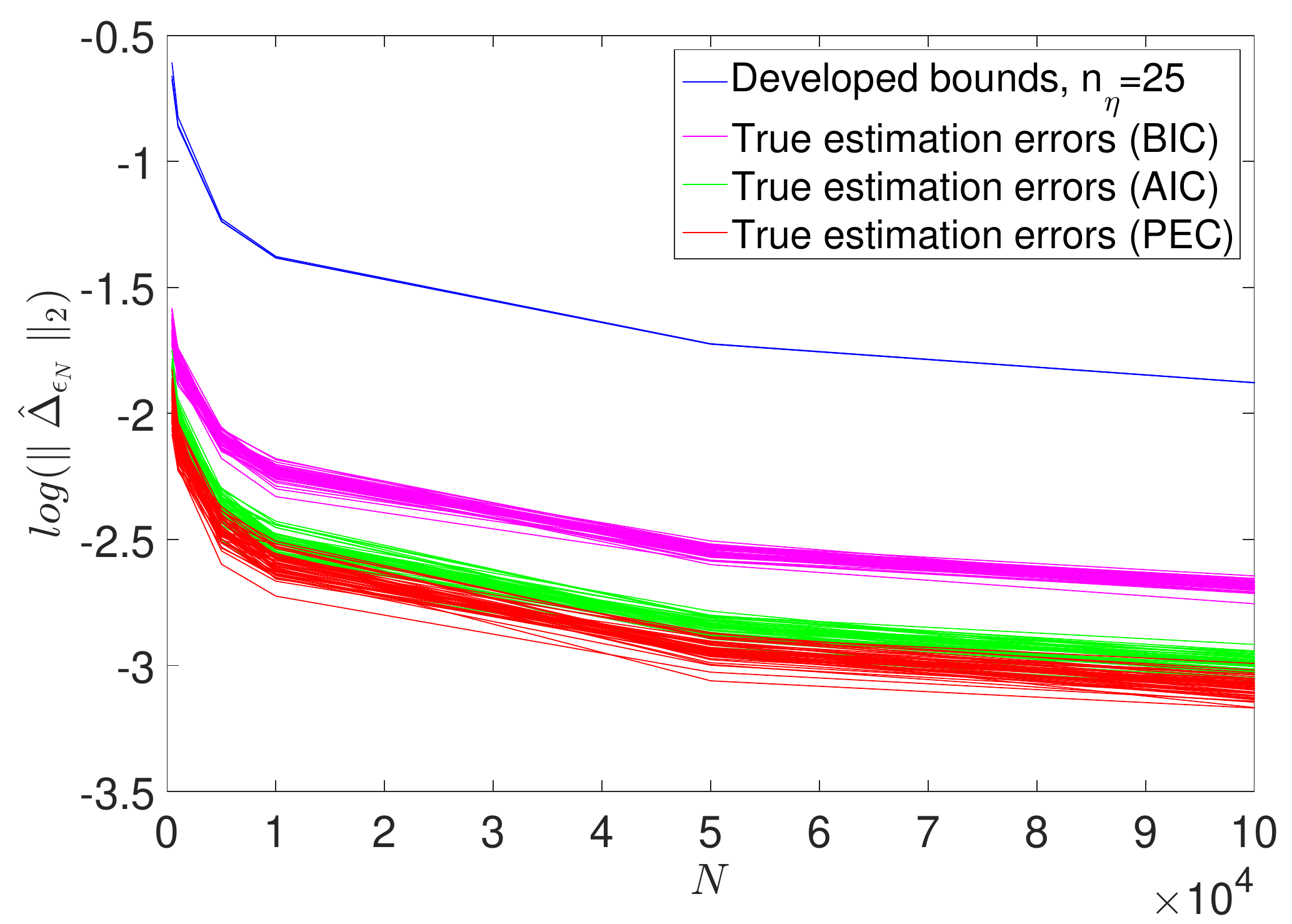}
\caption{Plot of the proposed bound and the true estimation errors corresponding to different choices of $\epsilon_N$, for a SNR=30dB and a white Gaussian input signal (magenta – BIC, green – AIC and red – PEC).} 
\label{fig_wn30dB_cb} 
\end{center}
\end{figure}

In addition, we plot another graph, shown in Fig. \ref{fig_wn30dB_cb}, to compare the proposed upper bound and the true estimation errors corresponding to different value of $\epsilon_N$, i.e. $\dfrac{n}{N}$ (PEC), $\dfrac{2n}{N}$ (AIC) and $\dfrac{log(N)n}{N}$ (BIC). The blue lines are the upper bounds developed in Theorem 5.1, which correspond to $n_\eta=25$, with the three different values of $\epsilon_N$. The magenta (BIC), green (AIC) and red (PEC) lines are the true estimation errors from 50 estimates (for each value of $\epsilon_N$) using the SPARSEVA framework. We can see that the plot again confirms the validity of the proposed upper bound for all choices of $\epsilon_N$. Note that the upper bound is not extremely tight, it is quite conservative, however, it is the price to usually pay for finite sample bounds with a general SPARSEVA setting, i.e. the regularized parameter $\epsilon_N$ can be any positive value. When $\epsilon_N$ is larger, the upper bound will be closer to the true estimate error.

%In addition, we also plot another graph, shown in Fig. , to compare the developed bound and the true estimate error corresponding to different values of $\epsilon_N$, e.g. $n/N$ (PEC), $2n/N$ (AIC) and $\text{log}(N)n/N$ (BIC).

\subsection{Coloured Noise Input}

In this section, a random discrete time system with a random model order between 1 and 10 is generated using the command \textit{drss} from Matlab. The system has poles with magnitude less than 0.9. White noise is added to the system output with different levels of SNR, e.g. 30dB, 20dB and 10dB. For each noise level, 50 different input excitation signals and output noise realizations are generated. For each set of input and output data, the system parameters are estimated using different sample sizes, i.e. $N=[450, 1000, 5000, 10000, 50000]$. 

Here, the input signal is generated by filtering a zero mean Gaussian white noise with unit variance through the filter,
$$F_u(q) = \dfrac{0.9798}{1-0.2q^{-1}}.$$ 
Due to this filtering, the covariance matrix of the regression matrix distribution will not be of a diagonal form. Note that this is a completely different scenario to that in Section 6.1. 

The FIR model structure is used here in order to construct the linear regression for the SPARSEVA problem (\ref{sparseva_ln}). The number of parameters $n$ of the FIR model is set to be 35. The regularization parameter, $\epsilon_N$, is chosen as $n/N$ \cite{Ha2015}.

We then compute the upper bound of $\Vert \hat{\theta}_{\epsilon_N}-\theta^* \Vert_2$ using (\ref{eq_tr_spbound}) with different values of $n_{\eta}$, i.e. $n_{\eta} = [10, 15, 25]$. The probability parameters $\alpha$ and $\beta$ are chosen to be $0.02$ and $0.001$ respectively. With this setting, the probability of the upper bound being correct is $(1-\alpha)(1-4n\beta)=0.84$. This upper bound will be compared with $\Vert \hat{\theta}_{\epsilon_N}-\theta^* \Vert_2$. %Note that here we plot both the upper bound and the true estimation errors in logarithmic scale.

\begin{figure}[h]
\begin{center}
\vspace{0.2cm}
\includegraphics[width=78mm,height=53mm]{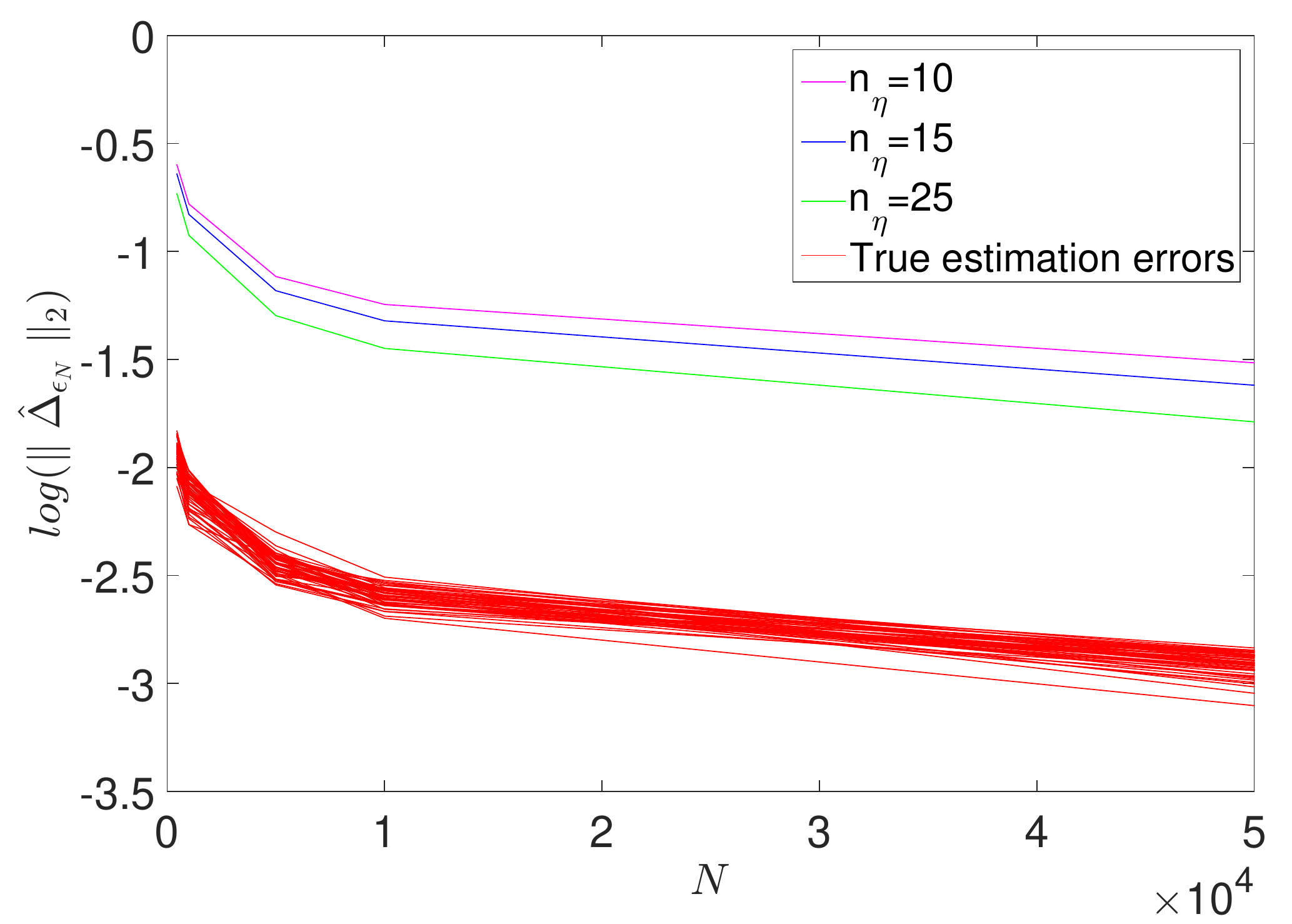}
\caption{Plot of the estimation error for SNR=30dB for a coloured input signal.} 
\label{fig_cl30dB} 
\end{center}
\end{figure}

\begin{figure}[h]
\begin{center}
\vspace{0.2cm}
\includegraphics[width=78mm,height=53mm]{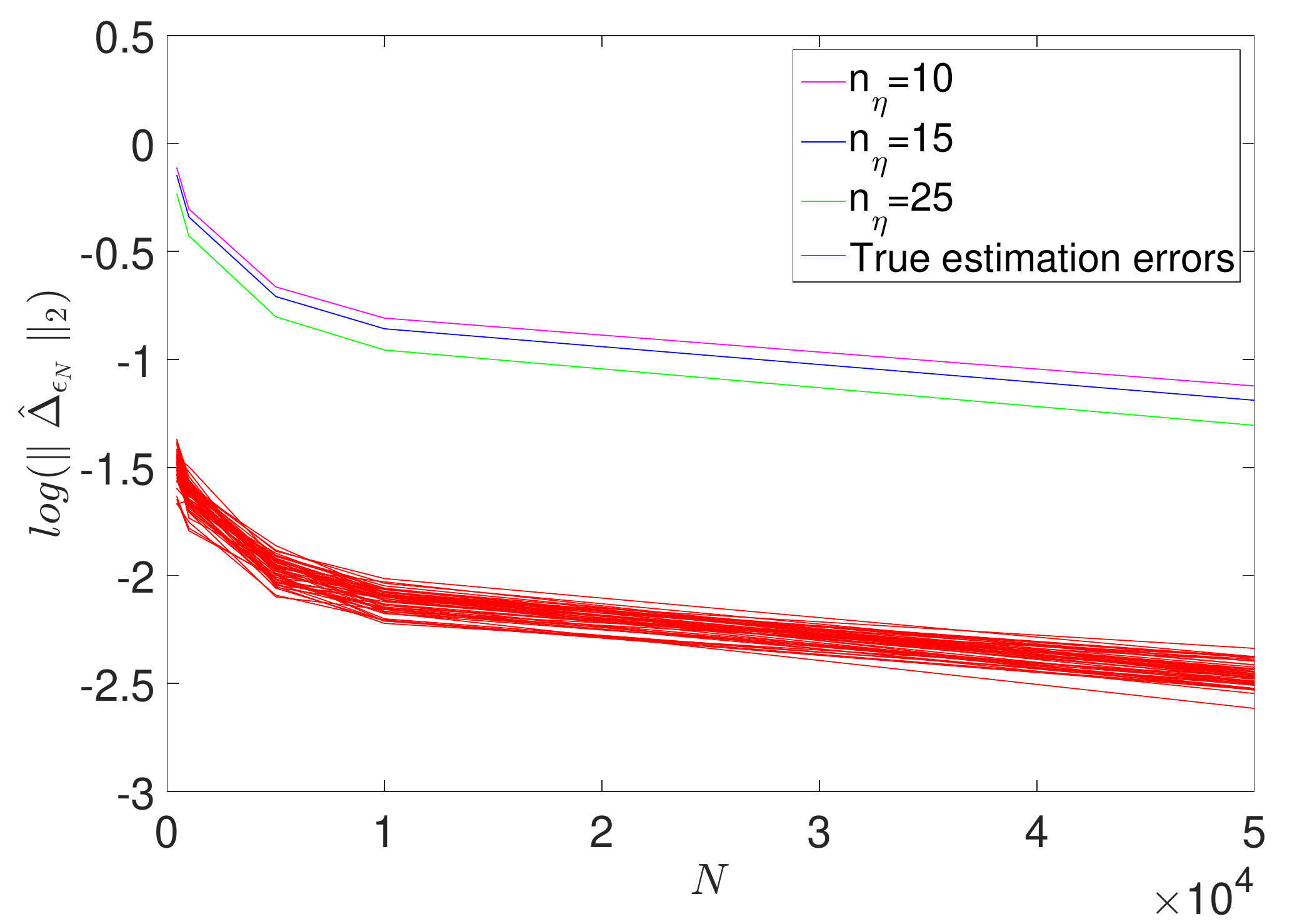}
\caption{Plot of the estimation error for SNR=20dB for a coloured input signal.} 
\label{fig_cl20dB} 
\end{center}
\end{figure}

\begin{figure}[h]
\begin{center}
\vspace{0.2cm}
\includegraphics[width=78mm,height=53mm]{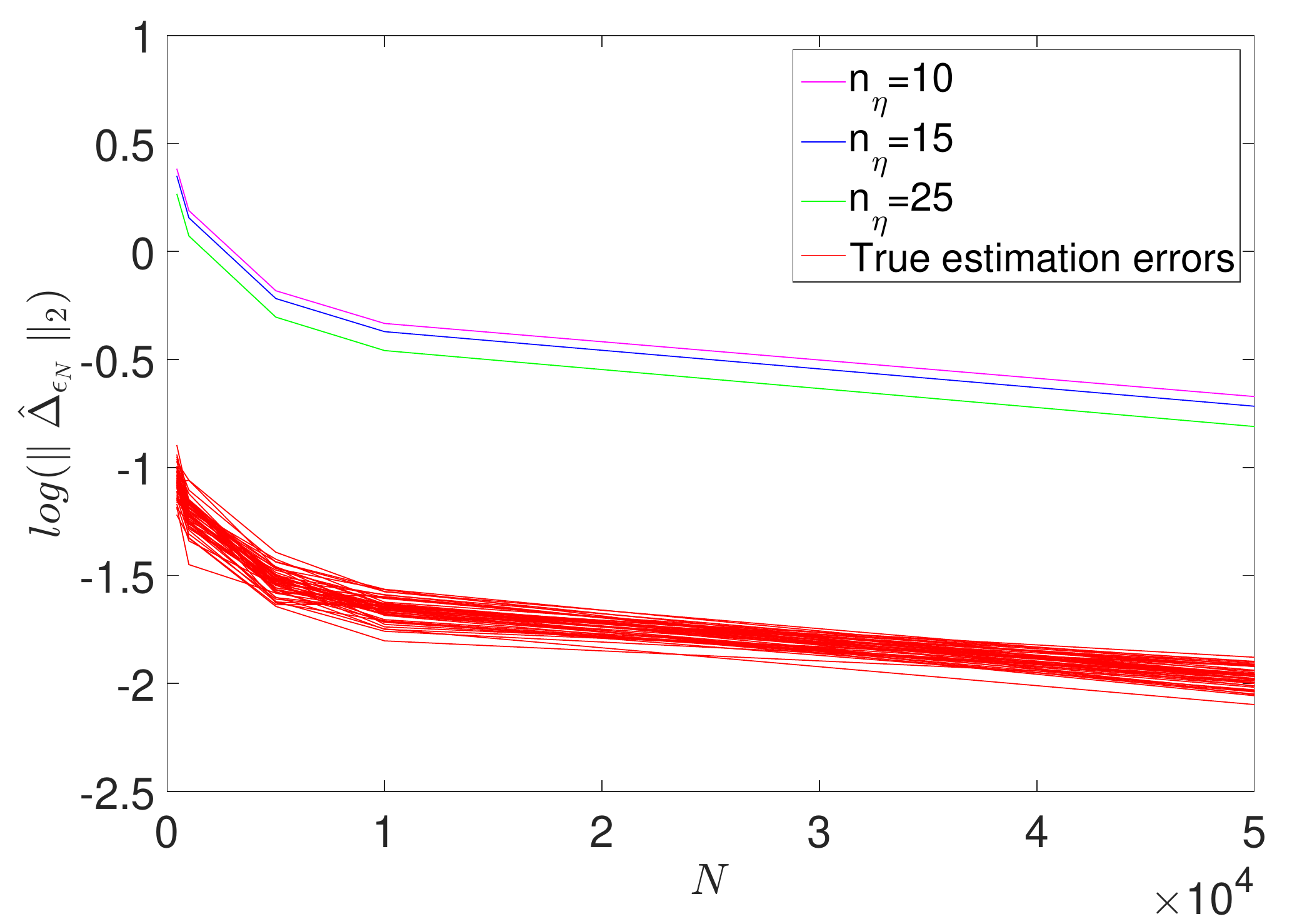}
\caption{Plot of the estimation error for SNR=10dB for a coloured input signal.} 
\label{fig_cl10dB} 
\end{center}
\end{figure}

Plots of the upper bound as stated in Theorem \ref{tr_spbound} and the true estimation error $\Vert \hat{\theta}_{\epsilon_N}-\theta^* \Vert_2$ are displayed in Figures \ref{fig_cl30dB} to \ref{fig_cl10dB}. In Figures \ref{fig_cl30dB} to \ref{fig_cl10dB}, the red lines are the true estimation errors from 50 estimates using the SPARSEVA framework. The magenta, blue and cyan lines are the upper bounds developed in Theorem \ref{tr_spbound}, which correspond to $n_\eta = [10,15,25]$ respectively. We can see that the plots confirmed the bound developed in Theorem \ref{tr_spbound} for all noise levels. When $N$ becomes large, the estimation error and the corresponding upper bound become smaller. When $N$ goes to infinity, the estimation error will tend to 0. Note that the bounds are slightly different for the chosen values of $n_\eta$, however, not significantly. As can be seen, the bounds are relatively insensitive to the choice of $n_\eta$.

%%%%%%%%%%%%%%%%%%%%%%%
\section{Conclusion}
%%%%%%%%%%%%%%%%%%%%%%%

The paper provides an upper bound on the SPARSEVA estimation error in the general case, for any choice of strongly convex cost function and decomposable norm. We also evaluate the bound for a specific scenario, i.e., a sparse regression estimate problem. Numerical results confirm the validity of the developed bound for different input signals with different output noise levels for different choices of the regularization parameters.

%
%\begin{ack}                               % Place acknowledgements
%Partially supported by the Roman Senate.  % here.
%\end{ack}

\bibliographystyle{plain}        % Include this if you use bibtex 
\bibliography{bibliography}           % and a bib file to produce the 
                                 % bibliography (preferred). The
                                 % correct style is generated by
                                 % Elsevier at the time of printing.

%\end{thebibliography}

\appendix

\section{Appendix}

\subsection{Background knowledge}

First, we cite a lemma directly from \cite{Negahban2012}, to enable the proof of Theorem \ref{tr_generalbound} to be constructed.
\begin{lemma} \label{lm_RInequality}
For any norm $\mathcal{R}$ that is decomposable with respect to $(\mathcal{M,\overline{\mathcal{M}}^{\perp}})$; and any vectors $\theta$, $\Delta$, we have
\begin{equation} \label{R_inequality}
\mathcal{R}(\theta+\Delta) - \mathcal{R}(\theta) \geq \mathcal{R}(\Delta_{\overline{\mathcal{M}}^{\perp}})-\mathcal{R}(\Delta_{\overline{\mathcal{M}}})-2\mathcal{R}(\theta_{\mathcal{M}^{\perp}}),
\end{equation}
\end{lemma} 
Recall that $\Delta_{\overline{\mathcal{M}}^{\perp}}$ is the Euclidean projection of $\Delta$ onto $\overline{\mathcal{M}}^{\perp}$(see Section~\ref{subsec:projection}), and similarly for the other terms in \eqref{R_inequality}.

\textit{\textbf{Proof.}} See the supplementary material of \cite{Negahban2012}. \hfill $\square$

We now quote the following lemma from \cite{Borwein2015} with modification to fit with the notation used in the SPARSEVA problem~(\ref{eq_c8sparseva}). This lemma helps us to find some important properties related to the SPARSEVA estimate $\hat{\theta}_{\lambda_N}$. Based on these properties, in the next section we can derive the upper bound on the estimation error. Note that for notational simplicity, we will denote $\mathcal{L}(\theta;Z_1^N)$ as $\mathcal{L}(\theta)$.

\begin{lemma} \label{lm_Lagrange}
Consider the convex optimization problem in (\ref{eq_c8sparseva}). Then the pair ($\hat{\theta}_{\epsilon_N}, \lambda_{\epsilon_N}$), with $\hat{\theta}_{\epsilon_N} \neq 0$, has the property that $\hat{\theta}_{\epsilon_N}$ is the solution of the problem (\ref{eq_c8sparseva}) and $\lambda_{\epsilon_N}$ is the Lagrange multiplier if and only if all of the following hold:
\begin{enumerate}
\item
$\lambda_{\epsilon_N} \in R^+$;
\item the function $\mathcal{R}(\theta) + \lambda\lbrace\mathcal{L}(\theta)-\mathcal{L}(\hat{\theta}_{NR})(1+\epsilon_N)\rbrace$ attains its minimum over $\mathbb{R}^n$ at $\hat{\theta}_{\epsilon_N}$; and
\item
$\mathcal{L}(\hat{\theta}_{\epsilon_N})-\mathcal{L}(\hat{\theta}_{NR})(1+\epsilon_N) = 0$.
\end{enumerate}
\end{lemma}
\textit{\textbf{Proof.}} See the proof of Theorem 2.2 in \cite{Borwein2015}. The third condition in the cited theorem is a complementary slackness condition, which reduces to condition (3) here if $\hat{\theta}_{\epsilon_N} \neq 0$ (cf. \cite[Lemma~II.2]{Rojas:14}).\hfill $\square$

\subsection{Proof of Theorem \ref{tr_generalbound}}

First we need to prove that there exists a Lagrange multiplier for the SPARSEVA problem (\ref{eq_c8sparseva}). We can assume without loss of generality that $\mathcal{L}(\hat{\theta}_{NR}) \neq 0$, since otherwise we can take $\lambda_{\epsilon_N} = 0$. According to \cite{Borwein2015}, the Lagrange multiplier for a convex optimization problem with constraint exists when the \textit{Slater condition} is satisfied. Specifically, for the SPARSEVA problem (\ref{eq_c8sparseva}), the Lagrange multiplier $\lambda_{\epsilon_N}$ exists when there exists a $\theta_1$ such that $\mathcal{L}(\theta_1) < \mathcal{L}(\hat{\theta}_{NR})(1+\epsilon_N)$. If $\epsilon_N > 0$, and $\mathcal{L}(\hat{\theta}_{NR}) \neq 0$, there always exists a parameter vector $\theta_1$ such that $\mathcal{L}(\theta_1) < \mathcal{L}(\hat{\theta}_{NR})(1+\epsilon_N)$ (just take $\theta_1 = \hat{\theta}_{NR}$). Therefore, there exists a Lagrange multiplier $\lambda_{\epsilon_N}$ for the SPARSEVA problem. 

Now that we have confirmed the existence of the Lagrange multiplier $\lambda_{\epsilon_N}$, consider the function $\mathcal{F}(\Delta)$ defined as follows,
\begin{equation} \label{eq_Fdelta_1}
\mathcal{F}(\Delta) = \mathcal{R}(\theta^*+\Delta) - \mathcal{R}(\theta^*) + \lambda_{\epsilon_N}\lbrace\mathcal{L}(\theta^*+\Delta)-\mathcal{L}(\theta^*)\rbrace.
\end{equation}
Using the strong convexity condition of $\mathcal{L}(\theta)$,
\begin{equation} \label{eq_sconvex}
\begin{aligned}
\mathcal{L}(\theta^*+\Delta)-\mathcal{L}(\theta^*)\geq \langle\nabla \mathcal{L}(\theta^*), \Delta\rangle + \kappa_{\mathcal{L}}\Vert \Delta \Vert^2_2.
\end{aligned}
\end{equation}
\CR{From \eqref{eq:dual_norm}, we have that}
%Applying the Cauchy-Schwartz inequality to $\mathcal{R}$ and its dual norm $\mathcal{R}^*$, we have, 
\begin{equation} \label{eq_CSIneqal}
\vert \langle\nabla \mathcal{L}(\theta^*), \Delta\rangle \vert \leq \mathcal{R}^*(\nabla\mathcal{L}(\theta^*))\mathcal{R}(\Delta).
\end{equation}
Next, combining the inequality (\ref{eq_CSIneqal}) and the triangle inequality, i.e., $\mathcal{R}(\Delta) \leq \mathcal{R}(\Delta_{\overline{\mathcal{M}}})+\mathcal{R}(\Delta_{\overline{\mathcal{M}}^\perp})$, we have,
\begin{equation} \nonumber
\begin{aligned}
\vert \langle\nabla \mathcal{L}(\theta^*), \Delta\rangle \vert & \leq \mathcal{R}^*(\nabla\mathcal{L}(\theta^*))\mathcal{R}(\Delta) \\
& \leq \mathcal{R}^*(\nabla\mathcal{L}(\theta^*))(\mathcal{R}(\Delta_{\overline{\mathcal{M}}})+\mathcal{R}(\Delta_{\overline{\mathcal{M}}^\perp})),
\end{aligned}
\end{equation}
therefore,
\begin{equation} \label{eq_dualnorm}
\begin{aligned}
\langle\nabla \mathcal{L}(\theta^*), \Delta\rangle \geq -\mathcal{R}^*(\nabla\mathcal{L}(\theta^*))(\mathcal{R}(\Delta_{\overline{\mathcal{M}}})+\mathcal{R}(\Delta_{\overline{\mathcal{M}}^\perp})).
\end{aligned}
\end{equation}
Now, combining (\ref{eq_Fdelta_1}), (\ref{eq_sconvex}), (\ref{eq_dualnorm}), and Lemma \ref{lm_RInequality}, 
\begin{align} \label{eq_Fdelta_gn}
& \mathcal{F}(\Delta) = \mathcal{R}(\theta^*+\Delta) - \mathcal{R}(\theta^*) + \lambda_{\epsilon_N}\lbrace\mathcal{L}(\theta^*+\Delta)-\mathcal{L}(\theta^*)\rbrace \nonumber \\
& \ \geq \mathcal{R}(\Delta_{\overline{\mathcal{M}}^{\perp}})-\mathcal{R}(\Delta_{\overline{\mathcal{M}}})-2\mathcal{R}(\theta^*_{\mathcal{M}^{\perp}}) \nonumber \\
& \ \ \ \ + \lambda_{\epsilon_N}\Big\lbrace -\mathcal{R}^*(\nabla\mathcal{L}(\theta^*)) (\mathcal{R}(\Delta_{\overline{\mathcal{M}}})+\mathcal{R}(\Delta_{\overline{\mathcal{M}}^\perp})) + \kappa_{\mathcal{L}}\Vert \Delta \Vert^2_2 \Big\rbrace \nonumber \\
& \ \geq \lbrace1-\lambda_{\epsilon_N}\mathcal{R}^*(\nabla\mathcal{L}(\theta^*))\rbrace\mathcal{R}(\Delta_{\overline{\mathcal{M}}^\perp}) \nonumber \\
& \ \ \ \ - \lbrace 1 + \lambda_{\epsilon_N}\mathcal{R}^*(\nabla\mathcal{L}(\theta^*))\rbrace\mathcal{R}(\Delta_{\overline{\mathcal{M}}}) + \kappa_{\mathcal{L}}\lambda_{\epsilon_N}\Vert \Delta \Vert^2_2 \nonumber \\
& \ \ \ \ -2\mathcal{R}(\theta^*_{\mathcal{M}^{\perp}}).
\end{align}

Notice that, when $\hat{\theta}_{\epsilon_N}$ is the estimate of the SPARSEVA problem, then property 2 in Lemma \ref{lm_Lagrange} states that the function $\mathcal{R}(\theta) + \lambda\lbrace\mathcal{L}(\theta)-\mathcal{L}(\hat{\theta}_{NR})(1+\epsilon_N)\rbrace$ attains its minimum over $\mathbb{R}^n$ at $\hat{\theta}_{\epsilon_N}$, which means,
\begin{multline*}
\forall \theta \in \mathcal{R}^n, \ \mathcal{R}(\hat{\theta}_{\epsilon_N}) + \lambda_{\epsilon_N}(\mathcal{L}(\hat{\theta}_{\epsilon_N})-\mathcal{L}(\hat{\theta}_{NR})(1+\epsilon_N)) \\
\leq \mathcal{R}(\theta) + \lambda_{\epsilon_N}(\mathcal{L}(\theta)-\mathcal{L}(\hat{\theta}_{NR})(1+\epsilon_N)).
\end{multline*}
Hence,
\begin{equation} \nonumber
\begin{aligned}
\forall \theta \in \mathcal{R}^n, \ & \mathcal{R}(\hat{\theta}_{\epsilon_N}) - \mathcal{R}(\theta) + \lambda_{\epsilon_N}\lbrace\mathcal{L}(\hat{\theta}_{\epsilon_N}) - \mathcal{L}(\theta)\rbrace \leq  0,
\end{aligned}
\end{equation}
or, taking $\theta = \theta^*$ and defining $\hat{\Delta}_{\epsilon_N} := \hat{\theta}_{\epsilon_N} - \theta^*$,
\begin{equation} \label{eq_Fdelta}
\mathcal{F}(\hat{\Delta}_{\epsilon_N}) \leq 0.
\end{equation}
Combining (\ref{eq_Fdelta_gn}) with (\ref{eq_Fdelta}), we then have,
\begin{equation} \label{eq_Fdelta_f}
\begin{aligned}
0 & \geq \kappa_{\mathcal{L}}\lambda_{\epsilon_N}\Vert \hat{\Delta}_{\epsilon_N} \Vert^2_2 + \lbrace1-\lambda_{\epsilon_N}\mathcal{R}^*(\nabla\mathcal{L}(\theta^*))\rbrace\mathcal{R}(\hat{\Delta}_{\epsilon_N,\overline{\mathcal{M}}^\perp}) \\
& \ \ \ \ \ - \lbrace 1 + \lambda_{\epsilon_N}\mathcal{R}^*(\nabla\mathcal{L}(\theta^*))\rbrace\mathcal{R}(\hat{\Delta}_{\epsilon_N,\overline{\mathcal{M}}})-2\mathcal{R}(\theta^*_{\mathcal{M}^{\perp}}).
\end{aligned}
\end{equation}

Now we consider two cases. \vspace{0.4cm} \\
\textbf{\underline{Case 1}:} \textbf{$\ \lambda_{\epsilon_N} \leq 1/\mathcal{R}^*(\nabla\mathcal{L}(\theta^*))$}

From (\ref{eq_Fdelta_f}), we have,
\begin{equation} \label{eq_Fdelta_c1}
\begin{aligned}
0 & \geq \kappa_{\mathcal{L}}\lambda_{\epsilon_N}\Vert \hat{\Delta}_{\epsilon_N} \Vert^2_2 - \lbrace 1 + \lambda_{\epsilon_N}\mathcal{R}^*(\nabla\mathcal{L}(\theta^*))\rbrace\mathcal{R}(\hat{\Delta}_{\epsilon_N,\overline{\mathcal{M}}}) \\
& \ \ \ \ -2\mathcal{R}(\theta^*_{\mathcal{M}^{\perp}}).
\end{aligned}
\end{equation}
By the definition of subspace compatibility, $$\mathcal{R}(\hat{\Delta}_{\epsilon_N,\overline{\mathcal{M}}}) \leq \Psi(\overline{\mathcal{M}})\Vert \hat{\Delta}_{\epsilon_N,\overline{\mathcal{M}}} \Vert_2.$$ 
Now \CR{we also have that}%since $0 \in \overline{\mathcal{M}}$, we have,
\begin{equation} \nonumber
%\begin{aligned}
\Vert \hat{\Delta}_{\epsilon_N,\overline{\mathcal{M}}} \Vert_2
= \CR{\Vert \Pi_{\overline{\mathcal{M}}}(\hat{\Delta}_{\epsilon_N})
%- \Pi_{\overline{\mathcal{M}}}(0)
\Vert_2}
\leq%\ \Vert \hat{\Delta}_{\epsilon_N} - 0 \Vert_2 = 
\CR{\Vert \hat{\Delta}_{\epsilon_N} \Vert_2.}
%\end{aligned}
\end{equation}
Therefore,
\begin{equation} \label{eq_subsComp}
\mathcal{R}(\hat{\Delta}_{\epsilon_N,\overline{\mathcal{M}}}) \leq \Psi(\overline{\mathcal{M}})\Vert \hat{\Delta}_{\epsilon_N} \Vert_2.
\end{equation}
Substituting this into (\ref{eq_Fdelta_c1}) gives,
\begin{multline} \label{eq_Fdelta_c1_2}
0 \geq \kappa_{\mathcal{L}} \lambda_{\epsilon_N} \Vert \hat{\Delta}_{\epsilon_N}  \Vert^2_2 - \lbrace 1 + \lambda_{\epsilon_N}\mathcal{R}^*(\nabla\mathcal{L}(\theta^*))\rbrace\Psi(\overline{\mathcal{M}})\Vert \hat{\Delta}_{\epsilon_N}  \Vert_2 \\
-2\mathcal{R}(\theta^*_{\mathcal{M}^{\perp}}).
\end{multline}
Note that for a quadratic polynomial $f(x) = ax^2+bx+c$, with $a > 0$, if there exists $x \in \mathbb{R}^+$ that makes $f(x) \leq 0$, then such $x$ must satisfy
\begin{equation} \nonumber
x \leq \dfrac{-b+\sqrt{b^2-4ac}}{2a}.
\end{equation}
Since $(A+B)^2 \leq 2A^2 + 2B^2$ for all $A,B \in \mathbb{R}$, 
\begin{equation} \label{eq_quadIneq}
x^2 \leq 2\left[ \dfrac{b^2}{4a^2}+\dfrac{b^2-4ac}{4a^2} \right] = \dfrac{b^2-2ac}{a^2}.
\end{equation}
Applying this inequality to (\ref{eq_Fdelta_c1_2}), we have,
\begin{align}
\Vert \hat{\Delta}_{\epsilon_N} \Vert^2_2\ & \leq \ \dfrac{1}{\kappa^2_{\mathcal{L}}\lambda^2_{\epsilon_N}}\lbrace 1 + \lambda_{\epsilon_N}\mathcal{R}^*(\nabla\mathcal{L}(\theta^*))\rbrace^2 \Psi^2(\overline{\mathcal{M}}) \nonumber \\
& \ \ \ \ + \dfrac{4}{\kappa_{\mathcal{L}}\lambda_{\epsilon_N}}\mathcal{R}(\theta^*_{\mathcal{M}^{\perp}}) \nonumber \\
%& \leq \ \dfrac{1}{\kappa^2_{\mathcal{L}}}\lbrace \dfrac{1}{\lambda_{\epsilon_N}} + \mathcal{R}^*(\nabla\mathcal{L}(\theta^*))\rbrace^2 \Psi^2(\overline{\mathcal{M}}) \\
%& \ \ \ \ + \dfrac{4}{\kappa_{\mathcal{L}}\lambda_{\epsilon_N}}\mathcal{R}(\theta^*_{\mathcal{M}^{\perp}}) \\
& \leq \ \dfrac{4}{\kappa^2_{\mathcal{L}}\lambda^2_{\epsilon_N}}\Psi^2(\overline{\mathcal{M}}) + \dfrac{4}{\kappa_{\mathcal{L}}\lambda_{\epsilon_N}}\mathcal{R}(\theta^*_{\mathcal{M}^{\perp}}).
\end{align}
\vspace{0.6cm}
\textbf{\underline{Case 2}:} $\lambda_{\epsilon_N} > 1/\mathcal{R}^*(\nabla\mathcal{L}(\theta^*))$ \\
Using a similar analysis as in Case 1,
\begin{equation} \label{eq_tempc2}
\mathcal{R}(\hat{\Delta}_{\epsilon_N,\overline{\mathcal{M}}^\perp}) \leq \Psi(\overline{\mathcal{M}}^\perp)\Vert \hat{\Delta}_{\epsilon_N} \Vert_2.
\end{equation}
Substituting (\ref{eq_tempc2}) and (\ref{eq_subsComp}) into (\ref{eq_Fdelta_f}), we obtain,
\begin{equation} \label{eq_Fdelta_c2}
\begin{aligned}
0 & \geq \kappa_{\mathcal{L}}\lambda_{\epsilon_N}\Vert \hat{\Delta}_{\epsilon_N} \Vert^2_2 + \lbrace1-\lambda_{\epsilon_N}\mathcal{R}^*(\nabla\mathcal{L}(\theta^*))\rbrace\Psi(\overline{\mathcal{M}}^\perp)\Vert \hat{\Delta}_{\epsilon_N} \Vert_2 \\
& \ \ \ - \lbrace 1 + \lambda_{\epsilon_N}\mathcal{R}^*(\nabla\mathcal{L}(\theta^*))\rbrace\Psi(\overline{\mathcal{M}})\Vert \hat{\Delta}_{\epsilon_N} \Vert_2 - 2\mathcal{R}(\theta^*_{\mathcal{M}^{\perp}}).
\end{aligned}
\end{equation}
Now using the inequality (\ref{eq_quadIneq}), yields,
\begin{equation} \label{eq_quadIneq2}
\begin{aligned}
\Vert \hat{\Delta}_{\epsilon_N} \Vert^2_2 & \leq \dfrac{1}{\kappa^2_{\mathcal{L}}} \Big(\left\lbrace\dfrac{1}{\lambda_{\epsilon_N}}-\mathcal{R}^*(\nabla\mathcal{L}(\theta^*))\right\rbrace\Psi(\overline{\mathcal{M}}^\perp) \\
& \ \ \ \ \ \ \ \ \ \ \ - \left\lbrace \dfrac{1}{\lambda_{\epsilon_N}} + \mathcal{R}^*(\nabla\mathcal{L}(\theta^*))\right\rbrace\Psi(\overline{\mathcal{M}}) \Big)^2 \\
& \ \ \ \ \ +\dfrac{4}{\kappa_{\mathcal{L}}\lambda_{\epsilon_N}}\mathcal{R}(\theta^*_{\mathcal{M}^{\perp}}).
\end{aligned}
\end{equation}
Applying the inequality $(A+B)^2 \leq 2A^2 + 2B^2$ to the first term in (\ref{eq_quadIneq2}) gives,
\begin{equation} \nonumber
\begin{aligned}
\Vert \hat{\Delta}_{\epsilon_N} \Vert^2_2 & \leq \dfrac{2}{\kappa^2_{\mathcal{L}}} \left\lbrace\dfrac{1}{\lambda_{\epsilon_N}}-\mathcal{R}^*(\nabla\mathcal{L}(\theta^*))\right\rbrace^2\Psi^2(\overline{\mathcal{M}}^\perp) \\
&\ \ \ \ + \dfrac{2}{\kappa^2_{\mathcal{L}}} \left\lbrace\dfrac{1}{\lambda_{\epsilon_N}}+\mathcal{R}^*(\nabla\mathcal{L}(\theta^*))\right\rbrace^2\Psi^2(\overline{\mathcal{M}}) \\
& \ \ \ \ +\dfrac{4}{\kappa_{\mathcal{L}}\lambda_{\epsilon_N}}\mathcal{R}(\theta^*_{\mathcal{M}^{\perp}}).
\end{aligned}
\end{equation}
Note that $0 < 1/\lambda_{\epsilon_N} < \mathcal{R}^*(\nabla\mathcal{L}(\theta^*))$, therefore,
\begin{equation} \label{eq_c8case2l}
\begin{aligned}
\left\lbrace \dfrac{1}{\lambda_{\epsilon_N}} - \mathcal{R}^*(\nabla\mathcal{L}(\theta^*))\right\rbrace^2 %
%& \leq \dfrac{1}{\lambda^2_{\epsilon_N}} + \lbrace\mathcal{R}^*(\nabla\mathcal{L}(\theta^*))\rbrace^2 \\
& \leq \lbrace\mathcal{R}^*(\nabla\mathcal{L}(\theta^*))\rbrace^2.
\end{aligned}
\end{equation}
We also have,
\begin{equation} \label{eq_c8case2l_2}
\begin{aligned}
\left\lbrace \dfrac{1}{\lambda_{\epsilon_N}} + \mathcal{R}^*(\nabla\mathcal{L}(\theta^*))\right\rbrace^2 \leq 4\lbrace\mathcal{R}^*(\nabla\mathcal{L}(\theta^*))\rbrace^2.
\end{aligned}
\end{equation}
Therefore, combining (\ref{eq_c8case2l}) and (\ref{eq_c8case2l_2}), 
\begin{equation} \nonumber
\begin{aligned}
\Vert \hat{\Delta}_{\epsilon_N} \Vert^2_2\ & \leq \ \dfrac{2}{\kappa^2_{\mathcal{L}}}\lbrace\mathcal{R}^*(\nabla\mathcal{L}(\theta^*))\rbrace^2 \Psi^2(\overline{\mathcal{M}}) \\
& \ \ +\dfrac{8}{\kappa^2_{\mathcal{L}}}\lbrace\mathcal{R}^*(\nabla\mathcal{L}(\theta^*))\rbrace^2 \Psi^2(\overline{\mathcal{M}}^\perp) \\
& \ \ + \dfrac{4}{\kappa_{\mathcal{L}}\lambda_{\epsilon_N}}\mathcal{R}(\theta^*_{\mathcal{M}^{\perp}}).
\end{aligned}
\end{equation} \hfill $\square$

\CR{
\subsection{Preliminary propositions for Theorem~\ref{tr_spbound}}

In this Appendix we present three propositions that assist in the development of the proof of Theorem~\ref{tr_spbound}:}

\begin{proposition} \label{pr_Lagrange}
Consider the optimization problem in (\ref{sparseva_ln}), and denote by $\lambda_{\epsilon_N}$ the corresponding Lagrange multiplier of its constraint. Then, if $\hat{\theta}_{\epsilon_N} \neq 0$, $\lambda_{\epsilon_N}$ can be computed as
\begin{equation} \label{eq_lagrange}
\lambda_{\epsilon_N}= \dfrac{1}{\Vert \nabla\mathcal{L}(\hat{\theta}_{\epsilon_N}) \Vert_\infty}.
\end{equation}
\end{proposition}
%\vspace{0.1cm}
\textit{\textbf{Proof.}} See Appendix A.4. \hfill $\square$

\begin{proposition} \label{pr_derThetastar}
Suppose Assumptions \ref{assu_X} and \ref{assu_e} hold. Then, with probability $1- n \beta$ ($0 \leq \beta \leq \CR{1/n}$), we have
\begin{equation} \nonumber
P\left( \left.\Vert \CR{\nabla \mathcal{L}(\theta^*)} \Vert_{\infty} \leq t \right| \Phi_N\right) \geq \left( 1 - 2\exp\left[ -\frac{N^2 t^2}{2 \sigma_e^2 \chi^2_{\beta}(\Sigma,I)} \right] \right)^n,
\end{equation}
where $s_{\text{max}}$ is the maximum element on the diagonal of the matrix $\Sigma$. In particular, \CR{choosing a specific value for $t$,}
\begin{equation} \label{eq_pr3spc}
\Vert \nabla \mathcal{L}(\theta^*) \Vert_{\infty} \leq \frac{\sqrt{2 \sigma_e^2 \chi^2_{\beta}(\Sigma,I)\ln(2/\beta)}}{N},
\end{equation}
with probability at least $1 - 2n \beta$ \CR{($0 \leq \beta \leq 1/2n$)}.
\end{proposition}
\textit{\textbf{Proof.}} See Appendix A.5. \hfill $\square$

\begin{proposition} \label{pr_derThetahat}
Suppose Assumptions \ref{assu_X} and \ref{assu_e} hold, then with probability at least $1- n \beta$ ($0 \leq \beta \leq \CR{1/n}$), we have
\begin{align*}
&P\left(\left. \Big\Vert \nabla\mathcal{L}(\hat{\theta}_{\epsilon_N}) \Big\Vert_{\infty} \leq t \right| e \right) \\
&\qquad\geq \left\lbrace 1 - 2\exp\left( -\frac{N^2 t^2}{2\sigma_e^2 s_{\text{max}} \chi^2_{\beta}(N-n) (1+\epsilon_N)} \right) \right\rbrace^n
\end{align*}
where $s_{\text{max}}$ is the maximum element on the diagonal of the matrix $\Sigma$. In particular, \CR{choosing a specific value for $t$,}
\begin{align} \label{eq_pr4spc}
\Big\Vert \nabla\mathcal{L}(\hat{\theta}_{\epsilon_N}) \Big\Vert_{\infty} \leq \frac{\sqrt{2\sigma_e^2 s_{\text{max}} \chi^2_{\beta}(N-n) (1+\epsilon_N)\ln(2/\beta)}}{N},
\end{align}
with probability at least $1 - 2 n \beta$ \CR{($0 \leq \beta \leq 1/2n$)}.
\end{proposition}
\textit{\textbf{Proof.}} See Appendix A.6. \hfill $\square$

\subsection{Proof of Proposition \ref{pr_Lagrange}}

Let us rewrite the SPARSEVA problem (\ref{sparseva_ln}),
\begin{equation} \label{sparseva_ln_sf}
\begin{aligned}
& \hat{\theta}_{\epsilon_N} \in \underset{\theta \in \mathbb{R}^n}{\text{arg min}}
& & \Vert\theta\Vert_1 \\
& \ \ \ \ \ \ \ \ \ \ \ \text{s.t.}
& & \mathcal{L}(\theta) - \mathcal{L}(\hat{\theta}_{NR})(1+\epsilon_N)\ \leq\ 0,
\end{aligned}
\end{equation}
in the Lagrangian form \eqref{eq_c8Mregularization} using Lemma \ref{lm_Lagrange}. The Lagrangian of the optimization problem (\ref{sparseva_ln_sf}) is,

\begin{equation} \label{eq_lagragrian}
\begin{aligned}
&g(\theta,\lambda) = \Vert \theta \Vert_1+ \lambda(\mathcal{L}(\theta) - \mathcal{L}(\hat{\theta}_{NR})(1+\epsilon_N)) \\
&= \Vert \theta \Vert_1 + \dfrac{\lambda}{2N}(\Vert Y_N-\Phi_N^T\theta \Vert^2_2 -\Vert Y_N-\Phi_N^T\hat{\theta}_{NR} \Vert^2_2(1+\epsilon_N)).
\end{aligned}
\end{equation}

The subdifferential of $g(\theta, \lambda)$ can be computed as
\begin{equation}
\frac{\partial g(\theta, \lambda)}{\partial \theta} = \text{v} - \dfrac{\lambda}{N} \Phi_N(Y_N-\Phi_N^T\theta),
\end{equation}
where $\text{v}=(v_1,\dots,v_m)^T$ is of the form
\begin{equation} \label{eq_vform}
\begin{aligned}
\begin{cases} v_i &= 1  \ \ \ \ \ \ \ \ \ \ \ \ \ \ \text{if} \ \ \theta_i>0 \\ v_i &= -1 \ \ \ \ \ \ \ \ \ \ \ \text{if} \ \ \theta_i<0 \\ v_i &\in[-1,1] \ \ \ \ \ \ \text{if} \ \ \theta_i=0.
\end{cases}
\end{aligned}
\end{equation}
Using property 2 of Lemma \ref{lm_Lagrange}, when $\hat{\theta}_{\epsilon_N}$ is a solution of the SPARSEVA problem (\ref{sparseva_ln}) and $\lambda_{\epsilon_N}$ is a Lagrange multiplier, we have,
\begin{equation}
\begin{aligned}
\textbf{0} & = \frac{\partial g(\theta, \lambda)}{\partial \theta} \Big\vert_{\theta=\hat{\theta}_{\epsilon_N}, \lambda=\lambda_{\epsilon_N}} \\
& =-\dfrac{\lambda_{\epsilon_N}}{N}\Phi_N(Y_N-\Phi_N^T\hat{\theta}_{\epsilon_N})+\text{v}_{\hat{\theta}_{\epsilon_N}},
\end{aligned}
\end{equation}
for some $\text{v}$ of the form in (\ref{eq_vform}). 
Note that when ${\hat{\theta}_{\epsilon_N}} \neq \textbf{0}$, $\Vert \text{v} \Vert_\infty=1$, which means that
$$ \lambda_{\epsilon_N} = \dfrac{N}{\Vert \Phi_N(Y_N-\Phi_N^T\hat{\theta}_{\epsilon_N}) \Vert_\infty}.$$
Since $\nabla\mathcal{L}(\hat{\theta}_{\epsilon_N})=\dfrac{1}{N}\Phi_N(Y_N-\Phi_N^T\hat{\theta}_{\epsilon_N})$, we can also write $\lambda_{\epsilon_N}$ as
$$ \lambda_{\epsilon_N}= \dfrac{1}{\Vert \nabla\mathcal{L}(\hat{\theta}_{\epsilon_N}) \Vert_\infty}.$$ \hfill $\square$

Note that this proof is similar to the one in \cite{Osborne2000}, where an expression was derived for the Lagrange multiplier in the traditional $l_1$ norm regularization problem (the LASSO). Here we have derived the Lagrange multiplier for the SPARSEVA problem as given in (\ref{sparseva_ln}). 

\subsection{Proof of Proposition \ref{pr_derThetastar}}

For the linear regression (\ref{eq_lr}) and the choice of $\mathcal{L}(\theta)$ in (\ref{eq_costFunc}), 
$$\nabla\mathcal{L}(\theta^*) = \dfrac{1}{N}\Phi_N(Y_N-\Phi_N^T\theta^*) = \dfrac{1}{N}\Phi_Ne.$$
Denote $R_j$ as the $j^{th}$ row of the matrix $\Phi_N$, then $ \nabla\mathcal{L}(\theta^*)$ can be computed as,
\begin{equation} \nonumber
\begin{aligned}
\nabla\mathcal{L}(\theta^*) = \dfrac{1}{N}\Phi_Ne = \dfrac{1}{N}
  \begin{bmatrix}
    R_1\\
    R_2\\
    \vdots \\
    R_n
  \end{bmatrix}e
  = \dfrac{1}{N}\begin{bmatrix}
    R_1e\\
    R_2e\\
    \vdots \\
    R_ne
  \end{bmatrix}
\end{aligned}
\end{equation}
consider the variable $Z = N^{-1} R_je$, using Assumption \ref{assu_e} on the disturbance noise $e$, $e \sim \mathcal{N}(0, \sigma_e^2)$, we have,
\begin{equation} \label{eq_DeltaLStarDist}
Z | e \sim \mathcal{N}\left(0, \dfrac{\sigma_e^2}{N^2}R_jR_j^T\right).
\end{equation}
Now in order to derive a bound for $\nabla\mathcal{L}(\theta^*)$, we first derive an upper bound for the variance $N^{-2} \sigma_e^2R_jR_j^T$ of the distribution in (\ref{eq_DeltaLStarDist}). Since $R_j \sim \mathcal{N}(0, \Sigma)$, we have,
\begin{equation} \nonumber
R_jR_j^T \sim \chi^2(\Sigma,\ I),
\end{equation}
where $\chi^2(\Sigma, I)$ is the generalized Chi squared with parameters $\Sigma$ and $I$. Hence, with probability $1-\beta,\ 0 \leq \beta \leq 1$, we have,
\begin{equation} \label{eq_XjSigmaXjT}
R_jR_j^T \leq \chi^2_{\beta}(\Sigma,\ I).
\end{equation}
Hence, the variance of the distribution of the variable $N^{-1} R_je$, is,
\begin{equation} \label{eq_varDeltaThetaStar_1}
\dfrac{\sigma_e^2}{N^2}R_jR_j^T \leq \dfrac{\sigma_e^2}{N^2}\chi^2_{\beta}(\Sigma,\ I),
\end{equation}
with probability $1-\beta$.

Note that from (\ref{eq_DeltaLStarDist}), for any $t > 0$, we have,
\begin{equation}
\begin{aligned}
P\left(\left. \Big \vert \dfrac{1}{N}R_je \Big \vert \leq t \right| \Phi_N \right) & = \int^t_{-t} f\left(x \left\vert 0, \dfrac{\sigma_e^2}{N^2}R_jR_j^T\right.\right)dx,
\end{aligned}
\end{equation}
where $f(x\vert 0,N^{-2} \sigma_e^2 R_jR_j^T)$ denotes the pdf of the Normal distribution $\mathcal{N}(0,N^{-2} \sigma_e^2 R_jR_j^T)$. This gives,
\begin{equation} \label{eq_probDeltaThetaStar}
P\left( \left.\Big\Vert \dfrac{\Phi_Ne}{N} \Big\Vert_{\infty} \leq t \right| \Phi_N\right) = \prod_{j=1}^{n} \left\lbrace \int^t_{-t} f\left(x \left\vert 0, \dfrac{\sigma_e^2}{N^2}R_jR_j^T\right.\right)dx \right\rbrace.
\end{equation}
This expression can be bounded from below using the standard result that $P(|\mathcal{N}(0,\sigma^2)| > t) \leq 2 \exp(-t^2/2\sigma^2)$ \cite[Eq.~(5.5)]{Vershynin-12}, to obtain
\begin{equation*}
P\left( \left.\Big\Vert \dfrac{\Phi_Ne}{N} \Big\Vert_{\infty} \leq t \right| \Phi_N\right) \geq \prod_{j=1}^{n} \left( 1 - 2\exp\left[ -\frac{N^2 t^2}{2 \sigma_e^2 R_jR_j^T} \right] \right).
\end{equation*}
The expression in parentheses on the right hand side is monotonically decreasing in $R_jR_j^T$, so using \eqref{eq_varDeltaThetaStar_1} gives
\begin{equation*}
P\left( \left.\Big\Vert \dfrac{\Phi_Ne}{N} \Big\Vert_{\infty} \leq t \right| \Phi_N\right) \geq \left( 1 - 2\exp\left[ -\frac{N^2 t^2}{2 \sigma_e^2 \chi^2_{\beta}(\Sigma,\ I)} \right] \right)^n,
\end{equation*}
which holds with probability\footnote{This bound follows because the events $A_j$ that \eqref{eq_varDeltaThetaStar_1} holds are not necessarily independent, but their joint probability can be bounded like $P(A_1 \cap \cdots \cap A_n) = 1 - P(A_1^C \cup \cdots \cup A_n^C) \geq 1 - P(A_1^C) - \cdots - P(A_n^C) = 1 - n\beta$.} at least $1- n \beta$.

In particular, taking $t = \sqrt{2 \sigma_e^2 \chi^2_{\beta}(\Sigma,\ I)\ln(2/\beta)} / N$ gives
\begin{align*}
P\left( \left.\Big\Vert \dfrac{\Phi_Ne}{N} \Big\Vert_{\infty} \leq \frac{\sqrt{2 \sigma_e^2 \chi^2_{\beta}(\Sigma,\ I)\ln(2/\beta)}}{N} \right| \Phi_N\right) %
&\geq (1 - \beta)^n \\
&\geq  1 - n \beta
\end{align*}
with probability at least $1- n \beta$, or equivalently,
\begin{align*}
\Big\Vert \nabla \mathcal{L}(\theta^*) \Big\Vert_{\infty} \leq \frac{\sqrt{2 \sigma_e^2 \chi^2_{\beta}(\Sigma,\ I)\ln(2/\beta)}}{N}
\end{align*}
with probability at least $1- 2 n \beta$. \hfill $\square$

\subsection{Proof of Proposition \ref{pr_derThetahat}}

When $\hat{\theta}_{\epsilon_N}$ is the solution of the problem in (\ref{sparseva_ln}), we have, 
\begin{equation} \label{eq_prl_1}
\nabla\mathcal{L}(\hat{\theta}_{\epsilon_N}) = \dfrac{1}{N}\Phi_N(Y_N-\Phi_N^T\hat{\theta}_{\epsilon_N}).
\end{equation}
Denote $e_{\epsilon_N} = Y_N-\Phi_N^T\hat{\theta}_{\epsilon_N}$, and $R_j$ as the $j^{th}$ row of the matrix $\Phi_N$, then (\ref{eq_prl_1}) becomes,
\begin{equation} \nonumber
\begin{aligned}
\nabla\mathcal{L}(\hat{\theta}_{\epsilon_N}) = \dfrac{1}{N}\Phi_N e_{\epsilon_N} = \dfrac{1}{N}
  \begin{bmatrix}
    R_1\\
    R_2\\
    \vdots \\
    R_n\\
  \end{bmatrix}e_{\epsilon_N}
  = \dfrac{1}{N}\begin{bmatrix}
    R_1e_{\epsilon_N}\\
    R_2e_{\epsilon_N}\\
    \vdots \\
    R_ne_{\epsilon_N}\\
  \end{bmatrix}.
  \end{aligned}
\end{equation}
From Assumption \ref{assu_X}, and using the same argument as in Proposition \ref{pr_derThetastar}, we have that each element of $R_j$ is distributed as $\mathcal{N}(0, \Sigma(j,j))$.

Consider the variable $Z=N^{-1}e_{\epsilon_N}^TR_j^T$, Since $R_j \sim \mathcal{N}(0,\ \Sigma)$,
\begin{equation} \nonumber
Z \sim \mathcal{N}\left(0, \dfrac{1}{N^2}e_{\epsilon_N}^T\Sigma e_{\epsilon_N}\right).
\end{equation}
Since $\Sigma$ is symmetric and positive definite matrix, hence using singular value decomposition, we can find a diagonal matrix $D$ that satisfies,
\begin{equation} \label{svd}
\Sigma = Q^TDQ,
\end{equation}
where $Q$ is the unitary matrix, i.e. $QQ^T=I$. Therefore, we have,
\begin{equation} \label{eq_varbound}
\dfrac{1}{N^2}e_{\epsilon_N}^T\Sigma e_{\epsilon_N} \leq \dfrac{s_{\text{max}}}{N^2} e_{\epsilon_N}^Te_{\epsilon_N},
\end{equation}
where $s_{\text{max}}$ is the maximum element on the diagonal of matrix $D$, i.e. maximum singular value of matrix $\Sigma$.
Note that, 
\begin{equation} \label{eq_eeNT2}
\begin{aligned}
e_{\epsilon_N}^Te_{\epsilon_N} &= (Y_N-\Phi_N^T\hat{\theta}_{\epsilon_N})^T(Y_N-\Phi_N^T\hat{\theta}_{\epsilon_N}) = 2N \mathcal{L}(\hat{\theta}_{\epsilon_N}) \\
&= 2N \mathcal{L}(\hat{\theta}_{NR})(1+\epsilon_N).
\end{aligned}
\end{equation}
From Section 4.4 in \cite{Soderstrom:89}, we have,
$$ \mathcal{L}(\hat{\theta}_{NR})|\Phi_N \sim \dfrac{\sigma_e^2}{2N} \chi^2(N-n),$$
which gives,
$$ \mathcal{L}(\hat{\theta}_{NR}) \leq \dfrac{\sigma_e^2}{2N} \chi^2_{\beta}(N-n),$$
with probability $1-\beta,\ 0 \leq \beta \leq 1$. Combining this inequality with (\ref{eq_varbound}) and (\ref{eq_eeNT2}) gives, with probability $1-\beta$,
$$ \dfrac{1}{N^2}e_{\epsilon_N}^T\Sigma e_{\epsilon_N} \leq \dfrac{\sigma_e^2}{N^2}s_{\text{max}} \chi^2_{\beta}(N-n) (1+\epsilon_N).$$
Hence,
\begin{multline} %\label{eq_probDeltaThetahat1}
P\left(\left. \left\vert \dfrac{R_je_{\epsilon_N}}{N} \right\vert \leq t \right| e \right) \\
\geq \int^t_{-t} f\left(x \left\vert 0, s_{\text{max}} \dfrac{\sigma_e^2}{N^2} \chi^2_{\beta}(N-n) (1+\epsilon_N)\right.\right)dx,
\end{multline}
with probability $1-\beta$. 
This means,
\begin{align} %\label{eq_probDeltaThetahat1}
&P\left(\left. \Big\Vert \dfrac{\Phi_Ne_{\epsilon_N}}{N} \Big\Vert_{\infty} \leq t \right| e \right) \\
&\qquad= \prod_{j=1}^n \left\lbrace \int^t_{-t} f\left(x \left\vert 0, \dfrac{\sigma_e^2}{N^2} s_{\text{max}} \chi^2_{\beta}(N-n) (1+\epsilon_N)\right.\right)dx \right\rbrace \nonumber \\
&\qquad\geq \left\lbrace 1 - 2\exp\left( -\frac{N^2 t^2}{2\sigma_e^2 s_{\text{max}} \chi^2_{\beta}(N-n) (1+\epsilon_N)} \right) \right\rbrace^n \nonumber
\end{align}
with probability at least $1- n \beta$, following the same reasoning as in the proof of Proposition~\ref{pr_derThetastar}. 
Therefore,
\begin{align*}
&P\left(\left. \Big\Vert \nabla\mathcal{L}(\hat{\theta}_{\epsilon_N}) \Big\Vert_{\infty} \leq t \right| e \right) \\
&\qquad\geq \left\lbrace 1 - 2\exp\left( -\frac{N^2 t^2}{2\sigma_e^2 s_{\text{max}} \chi^2_{\beta}(N-n) (1+\epsilon_N)} \right) \right\rbrace^n
\end{align*}
with probability at least $1 - n \beta$.

Taking $t = \sqrt{2\sigma_e^2 s_{\text{max}} \chi^2_{\beta}(N-n) (1+\epsilon_N)\ln(2/\beta)}/N$ gives
\begin{align*}
&P\left(\left. \Big\Vert \nabla\mathcal{L}(\hat{\theta}_{\epsilon_N}) \Big\Vert_{\infty} \leq \frac{\sqrt{2\sigma_e^2 s_{\text{max}} \chi^2_{\beta}(N-n) (1+\epsilon_N)\ln(2/\beta)}}{N} \right| e \right) \\
&\qquad\geq (1 - \beta)^n \geq 1 - n \beta,
\end{align*}
with probability at least $1 - n \beta$, or, equivalently,
\begin{align*}
\Big\Vert \nabla\mathcal{L}(\hat{\theta}_{\epsilon_N}) \Big\Vert_{\infty} \leq \frac{\sqrt{2\sigma_e^2 s_{\text{max}} \chi^2_{\beta}(N-n) (1+\epsilon_N)\ln(2/\beta)}}{N},
\end{align*}
with probability at least $1 - 2 n \beta$. \hfill $\square$

\end{document}